\def\filedate{20 November 2007}
\def\fileversion{2.1.9}

\documentclass{lms}
\usepackage{citeref}

\usepackage{lscape,graphicx}
\usepackage{graphicx}

\usepackage{amssymb}
\usepackage{amsmath}

\usepackage{tikz}
\usetikzlibrary{arrows}
\tikzstyle{block}=[draw opacity=0.7,line width=1.5cm]

\newcommand{\ndiv}{\hspace{-4pt}\not|\hspace{2pt}}

\newtheorem{theorem}{Theorem}[section] 

\newnumbered{assertion}{Assertion}    
\newnumbered{conjecture}{Conjecture}  
\newnumbered{definition}{Definition}
\newnumbered{hypothesis}{Hypothesis}
\newnumbered{remark}{Remark}
\newnumbered{note}{Note}
\newnumbered{observation}{Observation}
\newnumbered{problem}{Problem}
\newnumbered{question}{Question}
\newnumbered{algorithm}{Algorithm}
\newnumbered{example}{Example}
\newunnumbered{notation}{Notation} 

\title[On The Theory of the Collatz Number System]
 {On the Foundations of the Theory of a new Collatz Based Number System} 

\author{Michael A. Idowu}

\classno{00A71 (primary),  68R99, 03D75, 11B75 (secondary)} 

\extraline{Copyright \copyright Michael A. Idowu, 2014. All rights reserved.
}

\begin{document}
\bibliographystyle{plainnat}

\maketitle

\begin{abstract}
Set out here are some fundamental theories that may be regarded as newly discovered metamathematics of the odd integers in relation to the Collatz conjecture (also called the 3x+1 problem). Originally motivated by the requirement to invent a new optimised integer factorisation method, this foundational paper primarily focuses on the foundation, formalisation and presentation of a new theoretical framework (schema or blueprint) of a Collatz based number system. The proposed framework is based on metamathematical theories meticulously derived through iterative analyses and reverse engineering (i.e., by hand and mathematical computations) of many large subsets of integers. A collation of the fundamental results from these analytical attempts has led to the establishment of a completely deterministic model of a generalised Collatz based number system that is fundamentally and strangely associated with nonchaotic patterns. The proposed Collatz based number schema comprises of both visual and theoretical representations of many hidden patterns in Collatz sequences yet to be reported in literature. This novel theoretical approach may be viewed as a new method to contemporary Collatz conjecture research which may be connected to the proofs of many other mathematical theorems in number theory and discrete mathematics.
\end{abstract}



\section{Introduction}
Consider the set of natural numbers ${\bf \mathbb {N}^+} = \{1,2,3,4,...\} $, the Collatz mapping function $f$ maps each integer $n \in {\bf \mathbb N^+}$ to another positive integer also in ${\bf \mathbb  N^+}$ by the configuration
\begin{equation}
f(n)= \frac{3^{b(n)}n + b(n)}{2}
\end{equation}
where the binary flag function $b(n)=1$ is fixed if n is odd, otherwise $b(n)=0$ \cite{Ter76}. 

Note that this representation is equivalent to the Collatz mapping function 
\begin{equation}
f(n) = 
\begin{cases}
\frac{3n+1}{2} &\text{if  $n \equiv 1 $ (mod 2)}, i.e., ~[1]_2 \\
\frac{n}{2}  &\text{if $n \equiv 0$  (mod 2)}, i.e., ~ [0]_2 \\
\end{cases}
\end{equation}
often used in literature. The iterative execution of $f(n)$ always produces a sequence that terminates at 1 for any integer input $n$, according to the mathematician Lothar Collatz. This conjecture is known as the Collatz conjecture (CC) \cite{Lag85}, \cite{And98}, \cite{Van05} and was first proposed in 1937. It is commonly known as the \textit{3n+1} or \textit{3x+1} conjecture \cite{Lag85}, \cite{Marc96}, \cite{Gar81}, \cite{Sim99}, Ulam conjecture, Thwaites conjecture, the Kakutani's problem, the Hasse's algorithm, and the Syracuse problem \cite{Mad97}, \cite{Syr01}. This problem is easy to state but extremely difficult to prove. Many mathematicians have investigated and written articles about the conjetcure \cite{Ste77}, \cite{Bel06}, \cite{Bel98}, \cite{Sim05}, \cite{Sin10}. We recommend the work of Lagarias \cite{Lag85}, \cite{Lag11}, \cite{Lag12} for comprehensive and annotated bibliographies on the subject. 

The principal aim of this paper is to document, collate and present the results of a newly proposed number system formalised for investigating the truth of the CC.

\subsection{A new covering system and congruence classes modulo 18}
Here we seek to formulate a simple covering system for the odd integers such that every odd integer may be represented by exactly one of the residue classes $d_i \equiv r_i \mod18$ and each residue class is a finite (or infinite) collection (i.e., reordered set) of integers congruent to any member of the set $\{u_{i1}(\mathrm{mod}\ {2^13^2}),\ \ldots,\ u_{im}(\mathrm{mod}\ {2^m3^2})\}$. 

The covering system sought for the proposed Collatz number systems must be 1-cover, i.e., covers every integer exactly once. One of the objectives of this proposition is to develop a complete recategorisation of odd numbers based on both their residue classes and Collatz profiles (i.e., according to their $2^m $ divisibility properties after multiplying them by 3 and adding 1).

\subsection{Definitions}
Let $\sigma_{\infty}(n)$ represent the total stopping time of an integer $n$ under the Collatz system, where the symbol ``$\sigma_{\infty}$'' refers to ``number of iterations it takes to get to 1'' starting from the input $n$ - this is often referred to as the total stopping time of $n$.

Define the term \textit{Collatz profile} $f(d_i,m)$ as a representation of odd numbers $d_i$ congruent to $r_i \mod 18$ whose Collatz result $3d_i+1$ is divisible by $2^m$ such that the optimised result $\frac{3d_i+1}{2^m}$ is odd.  Hence, under the Collatz system the total stopping time relation $\sigma_{\infty}(d_i)=\sigma_{\infty}(\frac{3d_i+1}{2^m})+(m+1)$ holds, i.e., the total stopping time $\sigma_{\infty}(d_i)$ is $(m+1)$ iterates more than that of the next odd number in the Collatz sequence: $\sigma_{\infty}(\frac{3d+1}{2^m})=\sigma_{\infty}(d_i)-(m+1)$. 

For example, $f(13,3) = \frac{3(13)+1}{2^{\bf 3}} = 5$ implies that if the total stopping time of number {\bf 13} was $\sigma_{\infty}(13)$ under the Collatz system, then the total stopping time of number {\bf 5} \textit{would be} equal to $\sigma_{\infty}(5)=\sigma_{\infty}(13)-({\bf 3}+1)=\sigma_{\infty}(13)-4$, i.e., $\sigma_{\infty}(13) = 9$ and $\sigma_{\infty}(5) = 5$. Such simple inference (concept) may be used to prove the CC by \textit{reverse engineering}, i.e., the relative inference of total stopping times.

Define $r$ to be an odd number such that the following conditions are satisfied:
\begin{enumerate}
\item $1 \le r <18$;
\item $d_i=18V_{r_i, m}(n)+r_i$, $n \ge 0$;  and
\item $\frac{3d_i+1}{2^m}$ is odd.
 \end{enumerate}
If $V_{r, m}(n)$ is defined as the multiplicative set (function) that uniquely identifies numbers of same residue class attributable to similar Collatz-like pattern (profile), i.e., $V_{r, m}(n)$ is regarded as a parameter set dependent on actual values of $d_i$ and $r$ where $d_i$ is odd and congruent to $ r$ (mod 18), i.e., $r \in S$; where $S={1, 5, 3, 13, 17,15, 7, 11, 9}$. S is a reordered set of $r$. This set S establishes the connection between certain sets of integer numbers of one residue class in number system to other sets of other residue classes. The rationale behind this \textit{reordering} of the residue classes is of utmost importance (as demonstrated and explained in the next section). 

If we further stratify or reclassify every congruence class $d_i$ according to the result $ \frac{3d_i+1}{2^m}$, i.e., so that m is optimised and $ \frac{3d_i+1}{2^m}$ yields an odd number result, how does $d_i$ relate to the variable $m$, the exponent of $2^m$? Are there scientific or theoretical evidences to support the claim that the metatheory behind such relations may be completely deterministic? Ultimately, what is the schema for the proposed generalised Collatz based number system? These are the theoretical questions considered in this paper. 

From a theoretical point of view, the task of constructing a generalised schema of the Collatz based number system seems difficult, time consuming, and daunting. In practical terms, the difficulty associated with such tasks often seem more discouraging, because there is almost no guarantee of little or no recompense for the time invested. This is no longer the case. Here we present fundamental results that facilitate new theoretical perspectives on the Collatz conjecture and the generalised Collatz based number system. 

The next section (Section \ref{4np1}) introduces a modified (optimised) Collatz function and briefly gives a gentle introduction about a new discovery that inspires better understanding and new perspectives on residue classes modulo 18. 

The proposed metatheories of the Collatz based number system and its general proof are presented in section \ref{Reorder}, including the proposed schemata comprising of a map of generalised Collatz based number system and corresponding map of (distinct and fundamental) total stopping time functions.

\section{The fundamental relations between certain odd numbers} \label{4np1}
It is important to understand the fundamental relations between integers $d_i$ congruent to $r \mod 18$ and these relations define the inferred collatz properties, i.e., the divisibility quantities: $2^m |  (3d_i+1)$ analogous to the result $[p]_{54}$ where $1 \le p \le 54$ and p must be odd. In other words, how can new covering subsystems be formulated for the residue classes of odd integers, classifying each subsystem according to the prescribed Collatz properties? This question requires a thorough understanding of: a) the mechanisms of Collatz sequence transformations from odd to odd integers; b) the formulation of a new covering subsystems for the entire set of odd integers; and c) the determination of (total) stopping time for every odd integer for the production of an irrefutable proof of the Collatz conjecture. The last point, which is not addressed completely in this foundational paper, requires a thorough analysis of the propose theoretical schema of the Collatz based number system, i.e., a complete understanding of the covering system of all odd integers.

\subsection{The modified Collatz function}Introduce the modified Collatz function 
\begin{equation}\label{fdm}
f(d,m)= \frac{3^{b(d)}d + b(d)}{2^m} = \frac{3d+1}{2^m},
\end{equation}
where $m$ is the maximum exponent such that the odd number $d \in {\bf N^+}$ and the result $f(n,m)$ is an odd integer, i.e.,  $2^m \mid (3d+1)$ but $2^{m+1} \ndiv (3d+1)$ \cite{Ter76}.

\subsection{Residue classes modulo 18}
Given that $d$ is an odd number which belongs to only one of the following residue classes: 

\begin{equation}
d = 
\begin{cases}
d_1 &\text{if  $ d\equiv 1 $ (mod 18), e.g. 1, 349525, \dots}\\
d_2  &\text{if $ d \equiv 5$  (mod 18), e.g. 5}\\
\textcolor{red}{d_3  }&\text{if $ d \equiv 3$  (mod 18), e.g. \textcolor{red}{21}}\\
d_4  &\text{if $ d \equiv 13$  (mod 18), e.g. 85}\\
d_5  &\text{if $ d \equiv 17$  (mod 18), e.g. 341}\\
\textcolor{red}{d_6  }&\text{if $ d \equiv 15$  (mod 18), e.g. \textcolor{red}{1365}}\\
d_7  &\text{if $ d \equiv 7$  (mod 18), e.g. 5461}\\
d_8  &\text{if $ d \equiv 11$  (mod 18), e.g. 21845}\\
\textcolor{red}{d_9  }&\text{if $ d \equiv 9$  (mod 18), e.g. \textcolor{red}{87381}}\\
\hline
\textcolor{gray}{d_{10} }&\textcolor{gray}{\text{if $ d \equiv 1$  (mod 18), e.g. \textcolor{gray}{349525}}: \text{$ d_{10} = [1]_{18} = d_1$}} \\
\vdots
\end{cases}
\end{equation}
where $\{d_1 ~\cup ~d_2 ~\cup ~d_4 ~\cup ~d_5 \cup ~d_7 ~\cup ~d_8\}$ is the set of all odd numbers indivisible by 3, 
$[1]_{18} \in d_1$, 
$[5]_{18} \in d_2$, 
$[3]_{18} \in d_3$, 
$[13]_{18} \in d_4$, 
$[17]_{18} \in d_5$, 
$[15]_{18} \in d_6$, 
$[7]_{18} \in d_7$, 
$[11]_{18} \in d_8$, 
$[9]_{18} \in d_9$, and of course, 
$[0]_3 \in \{\{d_3~ \cup~d_6~ \cup~d_9\}$. This rearrangement enables one to capture cyclic recurrence relations between the odd integers using the relation $4d_i+1$, where $d_i \in d$ as demonstrated in the next subsection.

\subsection{Cyclic recurrence relations between odd numbers}
Let $r_i$ be a member of the finite set $S=\{1,5,3,13,17,15,7,11, 9\}$; $1 \le i \le 9$. This rearrangement of $[r_i]_{18}$ gives new insights into cyclic recurrence relations between the odd numbers:

\begin{equation}
d =
\begin{cases}
[r_i]_{18} &\text{if  $ 4([r_{i+8}]_{18})+1 \equiv r_{i} ~(mod ~ 18) $ \dots}\\
[r_{i+1}]_{18} &\text{if  $ 4([r_i]_{18})+1 \equiv r_{i+1} ~ (mod ~ 18) $ \dots}\\
\vdots \\
[r_{i+8}]_{18} &\text{if  $ 4([r_{i+7}]_{18})+1 \equiv r_{i+8} ~ (mod ~ 18)  $ }
\end{cases};
\end{equation}
 e.g. Let $ 19 \in d_i = [r_i]_{18} = [1]_{18} \equiv 1$ (mod 18), then  $ 77 = 4d_i+1 = [r_{i+1}]_{18} = [5]_{18} \equiv 5$ (mod 18),  $ 309 = 4^2d_i+5 = [r_{i+2}]_{18} = 3$ (mod 18),  $ 4^3d_i+21 = [r_{i+3}]_{18} $, $\dots$ , $ 1267029 = 4^8d_i+21845 =  [r_{i+8}]_{18}$, $ 5068117 = 4^9d_i+87381 =  [r_{i}]_{18}$, and so on. 
  
\section{The metatheory and properties of odd Integers }\label{Reorder}
We formulate and show how the arrays of $d_i$ that correspond to certain integers congruent to $r_i$ modulo 18 are defined, i.e., the $d_i=18V_{r_i, m}(n)+r_i$, $n \ge 0$ and $1 \le i \le 9$ that satisfy the required optimisation condition $\frac{3d_i+1}{2^m}=odd$. 

\subsection{New theories about certain odd numbers}
\begin{theorem}\label{mainthm}
Let the odd number transformation under the Collatz sequence system be defined as the function $f(d_i,m)= \frac{3d_i+1}{2^m} = d_{{i,m}_{next}}$ where $d_{{i,m}_{next}}$ represents the next odd number after the odd number $d_i$. The fundamental rules governing the relations between sets of $d_i$ and corresponding (i.e. mapped) sets of $V_{r_i,m}(n)$ are found to be completely deterministic (i.e. non-chaotic) on the following conditions enumerated:
\begin{enumerate}
\item
$d_i \equiv S(i) \mod 18$, i.e., $r_i = S(i)$;
\item
$d_i = 18V_{r_i,m}(n)+S(i)$; and
\item
the Collatz result $d_{{i,m}_{next}} = \frac{3(18~ V_{r_i,m}(n)+1)+S(i)}{2^m}$ must be odd.
\end{enumerate}
\end{theorem}

The proof of this major theorem requires enlisting all the possible $d_i$ and $V_{i,m}(n)$ values and demonstrating that the result $\frac{3(18V_{i,m}(n)+1)+S(i)}{2^m}$ is always odd. This theorem alone requires 162 explicit proofs - one for every statement. For the purpose of brevity, the trivialities involved in the explicit proofs are all avoided. 

\begin{proof}
When $ i =1, r_1 = 1$:

{\small
\begin{equation}
V_{1, m}(n) = 
\begin{cases}
\text{$V_{1,1} \in $ \{1,3,5,7,...\}} = & 2n + 1 \\
\text{$V_{1,2} \in $ \{0,4,8,12,...\}} = & 4n \\
\text{$V_{1,3} \in $ \{6,14,22,30,...\}} = & 8n + 6 \\
\text{$V_{1,4} \in $ \{2,18,34,50,...\}} = & 16n + 2 \\
\text{$V_{1,5} \in $ \{10,42,74,106,...\}} = & 32n + 10 \\
\text{$V_{1,6} \in $ \{58,122,186,250,...\}} = & 64n + 58 \\
\text{$V_{1,7} \in $ \{26,154,282,410,...\}} = & 128n + 26 \\
\text{$V_{1,8} \in $ \{90,346,602,858,...\}} = & 256n + 90 \\
\text{$V_{1,9} \in $ \{218,730,1242,1754,...\}} = & 512n + 218 \\
\text{$V_{1,10} \in $ \{474,1498,2522,3546,...\}} = & 1024n + 474 \\
\text{$V_{1,11} \in $ \{2010,4058,6106,8154,...\}} = & 2048n + 2010 \\
\text{$V_{1,12} \in $ \{986,5082,9178,13274,...\}} = & 4096n + 986 \\
\text{$V_{1,13} \in $ \{7130,15322,23514,31706,...\}} = & 8192n + 7130 \\
\text{$V_{1,14} \in $ \{11226,27610,43994,60378,...\}} = & 16384n + 11226 \\
\text{$V_{1,15} \in $ \{3034,35802,68570,101338,...\}} = & 32768n + 3034 \\
\text{$V_{1,16} \in $ \{52186,117722,183258,248794,...\}} = & 65536n + 52186 \\
\text{$V_{1,17} \in $ \{84954,216026,347098,478170,...\}} = & 131072n + 84954 \\
\text{$V_{1,18} \in $ \{150490,412634,674778,936922,...\}} = & 262144n + 150490 \\
\vdots
\end{cases}
\end{equation}
\begin{equation}
d_1 = 18V_{1, m}(n)+1 =
\begin{cases}
\text{$d_{1,1} \in $ \{19,55,91,127,...\}} = & 36n + 19 \\
\text{$d_{1,2} \in $ \{1,73,145,217,...\}} = & 72n + 1 \\
\text{$d_{1,3} \in $ \{109,253,397,541,...\}} = & 144n + 109 \\
\text{$d_{1,4} \in $ \{37,325,613,901,...\}} = & 288n + 37 \\
\text{$d_{1,5} \in $ \{181,757,1333,1909,...\}} = & 576n + 181 \\
\text{$d_{1,6} \in $ \{1045,2197,3349,4501,...\}} = & 1152n + 1045 \\
\text{$d_{1,7} \in $ \{469,2773,5077,7381,...\}} = & 2304n + 469 \\
\text{$d_{1,8} \in $ \{1621,6229,10837,15445,...\}} = & 4608n + 1621 \\
\text{$d_{1,9} \in $ \{3925,13141,22357,31573,...\}} = & 9216n + 3925 \\
\text{$d_{1,10} \in $ \{8533,26965,45397,63829,...\}} = & 18432n + 8533 \\
\text{$d_{1,11} \in $ \{36181,73045,109909,146773,...\}} = & 36864n + 36181 \\
\text{$d_{1,12} \in $ \{17749,91477,165205,238933,...\}} = & 73728n + 17749 \\
\text{$d_{1,13} \in $ \{128341,275797,423253,570709,...\}} = & 147456n + 128341 \\
\text{$d_{1,14} \in $ \{202069,496981,791893,1086805,...\}} = & 294912n + 202069 \\
\text{$d_{1,15} \in $ \{54613,644437,1234261,1824085,...\}} = & 589824n + 54613 \\
\text{$d_{1,16} \in $ \{939349,2118997,3298645,4478293,...\}} = & 1179648n + 939349 \\
\text{$d_{1,17} \in $ \{1529173,3888469,6247765,8607061,...\}} = & 2359296n + 1529173 \\
\text{$d_{1,18} \in $ \{2708821,7427413,12146005,16864597,...\}} = & 4718592n + 2708821 \\
\vdots
\end{cases}
\rightarrow \frac{3d_i+1}{2^m} = odd
\end{equation}
}

When $i =2, r_2 = S(2) = 5$ (note that the index of V below is 5):
{\small
\begin{equation}
V_{5, m}(n) = 
\begin{cases}
\text{$V_{5,1} \in $ \{1,3,5,7,...\}} = & 2n + 1 \\
\text{$V_{5,2} \in $ \{2,6,10,14,...\}} = & 4n + 2 \\
\text{$V_{5,3} \in $ \{4,12,20,28,...\}} = & 8n + 4 \\
\text{$V_{5,4} \in $ \{0,16,32,48,...\}} = & 16n \\
\text{$V_{5,5} \in $ \{24,56,88,120,...\}} = & 32n + 24 \\
\text{$V_{5,6} \in $ \{8,72,136,200,...\}} = & 64n + 8 \\
\text{$V_{5,7} \in $ \{40,168,296,424,...\}} = & 128n + 40 \\
\text{$V_{5,8} \in $ \{232,488,744,1000,...\}} = & 256n + 232 \\
\text{$V_{5,9} \in $ \{104,616,1128,1640,...\}} = & 512n + 104 \\
\text{$V_{5,10} \in $ \{360,1384,2408,3432,...\}} = & 1024n + 360 \\
\text{$V_{5,11} \in $ \{872,2920,4968,7016,...\}} = & 2048n + 872 \\
\text{$V_{5,12} \in $ \{1896,5992,10088,14184,...\}} = & 4096n + 1896 \\
\text{$V_{5,13} \in $ \{8040,16232,24424,32616,...\}} = & 8192n + 8040 \\
\text{$V_{5,14} \in $ \{3944,20328,36712,53096,...\}} = & 16384n + 3944 \\
\text{$V_{5,15} \in $ \{28520,61288,94056,126824,...\}} = & 32768n + 28520 \\
\text{$V_{5,16} \in $ \{44904,110440,175976,241512,...\}} = & 65536n + 44904 \\
\text{$V_{5,17} \in $ \{12136,143208,274280,405352,...\}} = & 131072n + 12136 \\
\text{$V_{5,18} \in $ \{208744,470888,733032,995176,...\}} = & 262144n + 208744 \\
\vdots
\end{cases}
\end{equation}
\begin{equation}
d_2 = 18V_{5, m}(n)+5 = 
\begin{cases}
\text{$d_{2,1} \in $ \{23,59,95,131,...\}} = & 36n + 23 \\
\text{$d_{2,2} \in $ \{41,113,185,257,...\}} = & 72n + 41 \\
\text{$d_{2,3} \in $ \{77,221,365,509,...\}} = & 144n + 77 \\
\text{$d_{2,4} \in $ \{5,293,581,869,...\}} = & 288n + 5 \\
\text{$d_{2,5} \in $ \{437,1013,1589,2165,...\}} = & 576n + 437 \\
\text{$d_{2,6} \in $ \{149,1301,2453,3605,...\}} = & 1152n + 149 \\
\text{$d_{2,7} \in $ \{725,3029,5333,7637,...\}} = & 2304n + 725 \\
\text{$d_{2,8} \in $ \{4181,8789,13397,18005,...\}} = & 4608n + 4181 \\
\text{$d_{2,9} \in $ \{1877,11093,20309,29525,...\}} = & 9216n + 1877 \\
\text{$d_{2,10} \in $ \{6485,24917,43349,61781,...\}} = & 18432n + 6485 \\
\text{$d_{2,11} \in $ \{15701,52565,89429,126293,...\}} = & 36864n + 15701 \\
\text{$d_{2,12} \in $ \{34133,107861,181589,255317,...\}} = & 73728n + 34133 \\
\text{$d_{2,13} \in $ \{144725,292181,439637,587093,...\}} = & 147456n + 144725 \\
\text{$d_{2,14} \in $ \{70997,365909,660821,955733,...\}} = & 294912n + 70997 \\
\text{$d_{2,15} \in $ \{513365,1103189,1693013,2282837,...\}} = & 589824n + 513365 \\
\text{$d_{2,16} \in $ \{808277,1987925,3167573,4347221,...\}} = & 1179648n + 808277 \\
\text{$d_{2,17} \in $ \{218453,2577749,4937045,7296341,...\}} = & 2359296n + 218453 \\
\text{$d_{2,18} \in $ \{3757397,8475989,13194581,17913173,...\}} = & 4718592n + 3757397 \\
\vdots
\end{cases}
\rightarrow \frac{3d_2+1}{2^m} = odd
\end{equation}
}

When $i =3, r_3 = S(3) = 3$:
{\small
\begin{equation}
V_{3, m}(n) = 
\begin{cases}
\text{$V_{3,1} \in $ \{0,2,4,6,...\}} = & 2n \\
\text{$V_{3,2} \in $ \{3,7,11,15,...\}} = & 4n + 3 \\
\text{$V_{3,3} \in $ \{5,13,21,29,...\}} = & 8n + 5 \\
\text{$V_{3,4} \in $ \{9,25,41,57,...\}} = & 16n + 9 \\
\text{$V_{3,5} \in $ \{17,49,81,113,...\}} = & 32n + 17 \\
\text{$V_{3,6} \in $ \{1,65,129,193,...\}} = & 64n + 1 \\
\text{$V_{3,7} \in $ \{97,225,353,481,...\}} = & 128n + 97 \\
\text{$V_{3,8} \in $ \{33,289,545,801,...\}} = & 256n + 33 \\
\text{$V_{3,9} \in $ \{161,673,1185,1697,...\}} = & 512n + 161 \\
\text{$V_{3,10} \in $ \{929,1953,2977,4001,...\}} = & 1024n + 929 \\
\text{$V_{3,11} \in $ \{417,2465,4513,6561,...\}} = & 2048n + 417 \\
\text{$V_{3,12} \in $ \{1441,5537,9633,13729,...\}} = & 4096n + 1441 \\
\text{$V_{3,13} \in $ \{3489,11681,19873,28065,...\}} = & 8192n + 3489 \\
\text{$V_{3,14} \in $ \{7585,23969,40353,56737,...\}} = & 16384n + 7585 \\
\text{$V_{3,15} \in $ \{32161,64929,97697,130465,...\}} = & 32768n + 32161 \\
\text{$V_{3,16} \in $ \{15777,81313,146849,212385,...\}} = & 65536n + 15777 \\
\text{$V_{3,17} \in $ \{114081,245153,376225,507297,...\}} = & 131072n + 114081 \\
\text{$V_{3,18} \in $ \{179617,441761,703905,966049,...\}} = & 262144n + 179617 \\
\vdots
\end{cases}
\end{equation}
\begin{equation}
d_3 = 18V_{3, m}(n)+3 =
\begin{cases}
\text{$d_{3,1} \in $ \{3,39,75,111,...\}} = & 36n + 3 \\
\text{$d_{3,2} \in $ \{57,129,201,273,...\}} = & 72n + 57 \\
\text{$d_{3,3} \in $ \{93,237,381,525,...\}} = & 144n + 93 \\
\text{$d_{3,4} \in $ \{165,453,741,1029,...\}} = & 288n + 165 \\
\text{$d_{3,5} \in $ \{309,885,1461,2037,...\}} = & 576n + 309 \\
\text{$d_{3,6} \in $ \{21,1173,2325,3477,...\}} = & 1152n + 21 \\
\text{$d_{3,7} \in $ \{1749,4053,6357,8661,...\}} = & 2304n + 1749 \\
\text{$d_{3,8} \in $ \{597,5205,9813,14421,...\}} = & 4608n + 597 \\
\text{$d_{3,9} \in $ \{2901,12117,21333,30549,...\}} = & 9216n + 2901 \\
\text{$d_{3,10} \in $ \{16725,35157,53589,72021,...\}} = & 18432n + 16725 \\
\text{$d_{3,11} \in $ \{7509,44373,81237,118101,...\}} = & 36864n + 7509 \\
\text{$d_{3,12} \in $ \{25941,99669,173397,247125,...\}} = & 73728n + 25941 \\
\text{$d_{3,13} \in $ \{62805,210261,357717,505173,...\}} = & 147456n + 62805 \\
\text{$d_{3,14} \in $ \{136533,431445,726357,1021269,...\}} = & 294912n + 136533 \\
\text{$d_{3,15} \in $ \{578901,1168725,1758549,2348373,...\}} = & 589824n + 578901 \\
\text{$d_{3,16} \in $ \{283989,1463637,2643285,3822933,...\}} = & 1179648n + 283989 \\
\text{$d_{3,17} \in $ \{2053461,4412757,6772053,9131349,...\}} = & 2359296n + 2053461 \\
\text{$d_{3,18} \in $ \{3233109,7951701,12670293,17388885,...\}} = & 4718592n + 3233109 \\
\vdots
\end{cases}
\rightarrow \frac{3d_3+1}{2^m} = odd
\end{equation}
}

When $ i =4, r_4 = S(4) = 13$ (note that the index of V below is 13):
{\small
\begin{equation}
V_{13, m}(n) = 
\begin{cases}
\text{$V_{13,1} \in $ \{1,3,5,7,...\}} = & 2n + 1 \\
\text{$V_{13,2} \in $ \{2,6,10,14,...\}} = & 4n + 2 \\
\text{$V_{13,3} \in $ \{0,8,16,24,...\}} = & 8n \\
\text{$V_{13,4} \in $ \{12,28,44,60,...\}} = & 16n + 12 \\
\text{$V_{13,5} \in $ \{20,52,84,116,...\}} = & 32n + 20 \\
\text{$V_{13,6} \in $ \{36,100,164,228,...\}} = & 64n + 36 \\
\text{$V_{13,7} \in $ \{68,196,324,452,...\}} = & 128n + 68 \\
\text{$V_{13,8} \in $ \{4,260,516,772,...\}} = & 256n + 4 \\
\text{$V_{13,9} \in $ \{388,900,1412,1924,...\}} = & 512n + 388 \\
\text{$V_{13,10} \in $ \{132,1156,2180,3204,...\}} = & 1024n + 132 \\
\text{$V_{13,11} \in $ \{644,2692,4740,6788,...\}} = & 2048n + 644 \\
\text{$V_{13,12} \in $ \{3716,7812,11908,16004,...\}} = & 4096n + 3716 \\
\text{$V_{13,13} \in $ \{1668,9860,18052,26244,...\}} = & 8192n + 1668 \\
\text{$V_{13,14} \in $ \{5764,22148,38532,54916,...\}} = & 16384n + 5764 \\
\text{$V_{13,15} \in $ \{13956,46724,79492,112260,...\}} = & 32768n + 13956 \\
\text{$V_{13,16} \in $ \{30340,95876,161412,226948,...\}} = & 65536n + 30340 \\
\text{$V_{13,17} \in $ \{128644,259716,390788,521860,...\}} = & 131072n + 128644 \\
\text{$V_{13,18} \in $ \{63108,325252,587396,849540,...\}} = & 262144n + 63108 \\
\vdots
\end{cases}
\end{equation}
\begin{equation}
d_4 = 18V_{13, m}(n)+13 =
\begin{cases}
\text{$d_{4,1} \in $ \{31,67,103,139,...\}} = & 36n + 31 \\
\text{$d_{4,2} \in $ \{49,121,193,265,...\}} = & 72n + 49 \\
\text{$d_{4,3} \in $ \{13,157,301,445,...\}} = & 144n + 13 \\
\text{$d_{4,4} \in $ \{229,517,805,1093,...\}} = & 288n + 229 \\
\text{$d_{4,5} \in $ \{373,949,1525,2101,...\}} = & 576n + 373 \\
\text{$d_{4,6} \in $ \{661,1813,2965,4117,...\}} = & 1152n + 661 \\
\text{$d_{4,7} \in $ \{1237,3541,5845,8149,...\}} = & 2304n + 1237 \\
\text{$d_{4,8} \in $ \{85,4693,9301,13909,...\}} = & 4608n + 85 \\
\text{$d_{4,9} \in $ \{6997,16213,25429,34645,...\}} = & 9216n + 6997 \\
\text{$d_{4,10} \in $ \{2389,20821,39253,57685,...\}} = & 18432n + 2389 \\
\text{$d_{4,11} \in $ \{11605,48469,85333,122197,...\}} = & 36864n + 11605 \\
\text{$d_{4,12} \in $ \{66901,140629,214357,288085,...\}} = & 73728n + 66901 \\
\text{$d_{4,13} \in $ \{30037,177493,324949,472405,...\}} = & 147456n + 30037 \\
\text{$d_{4,14} \in $ \{103765,398677,693589,988501,...\}} = & 294912n + 103765 \\
\text{$d_{4,15} \in $ \{251221,841045,1430869,2020693,...\}} = & 589824n + 251221 \\
\text{$d_{4,16} \in $ \{546133,1725781,2905429,4085077,...\}} = & 1179648n + 546133 \\
\text{$d_{4,17} \in $ \{2315605,4674901,7034197,9393493,...\}} = & 2359296n + 2315605 \\
\text{$d_{4,18} \in $ \{1135957,5854549,10573141,15291733,...\}} = & 4718592n + 1135957 \\
\vdots
\end{cases}
\rightarrow \frac{3d_4+1}{2^m} = odd
\end{equation}
}

When $i =5, r_5 = S(5) = 17$ (note that the index of V below is 17):
{\small
\begin{equation}
V_{17, m}(n) = 
\begin{cases}
\text{$V_{17,1} \in $ \{1,3,5,7,...\}} = & 2n + 1 \\
\text{$V_{17,2} \in $ \{0,4,8,12,...\}} = & 4n \\
\text{$V_{17,3} \in $ \{6,14,22,30,...\}} = & 8n + 6 \\
\text{$V_{17,4} \in $ \{10,26,42,58,...\}} = & 16n + 10 \\
\text{$V_{17,5} \in $ \{2,34,66,98,...\}} = & 32n + 2 \\
\text{$V_{17,6} \in $ \{50,114,178,242,...\}} = & 64n + 50 \\
\text{$V_{17,7} \in $ \{82,210,338,466,...\}} = & 128n + 82 \\
\text{$V_{17,8} \in $ \{146,402,658,914,...\}} = & 256n + 146 \\
\text{$V_{17,9} \in $ \{274,786,1298,1810,...\}} = & 512n + 274 \\
\text{$V_{17,10} \in $ \{18,1042,2066,3090,...\}} = & 1024n + 18 \\
\text{$V_{17,11} \in $ \{1554,3602,5650,7698,...\}} = & 2048n + 1554 \\
\text{$V_{17,12} \in $ \{530,4626,8722,12818,...\}} = & 4096n + 530 \\
\text{$V_{17,13} \in $ \{2578,10770,18962,27154,...\}} = & 8192n + 2578 \\
\text{$V_{17,14} \in $ \{14866,31250,47634,64018,...\}} = & 16384n + 14866 \\
\text{$V_{17,15} \in $ \{6674,39442,72210,104978,...\}} = & 32768n + 6674 \\
\text{$V_{17,16} \in $ \{23058,88594,154130,219666,...\}} = & 65536n + 23058 \\
\text{$V_{17,17} \in $ \{55826,186898,317970,449042,...\}} = & 131072n + 55826 \\
\text{$V_{17,18} \in $ \{121362,383506,645650,907794,...\}} = & 262144n + 121362 \\
\vdots
\end{cases}
\end{equation}
\begin{equation}
d_5 = 18V_{17, m}(n) +17 =
\begin{cases}
\text{$d_{5,1} \in $ \{35,71,107,143,...\}} = & 36n + 35 \\
\text{$d_{5,2} \in $ \{17,89,161,233,...\}} = & 72n + 17 \\
\text{$d_{5,3} \in $ \{125,269,413,557,...\}} = & 144n + 125 \\
\text{$d_{5,4} \in $ \{197,485,773,1061,...\}} = & 288n + 197 \\
\text{$d_{5,5} \in $ \{53,629,1205,1781,...\}} = & 576n + 53 \\
\text{$d_{5,6} \in $ \{917,2069,3221,4373,...\}} = & 1152n + 917 \\
\text{$d_{5,7} \in $ \{1493,3797,6101,8405,...\}} = & 2304n + 1493 \\
\text{$d_{5,8} \in $ \{2645,7253,11861,16469,...\}} = & 4608n + 2645 \\
\text{$d_{5,9} \in $ \{4949,14165,23381,32597,...\}} = & 9216n + 4949 \\
\text{$d_{5,10} \in $ \{341,18773,37205,55637,...\}} = & 18432n + 341 \\
\text{$d_{5,11} \in $ \{27989,64853,101717,138581,...\}} = & 36864n + 27989 \\
\text{$d_{5,12} \in $ \{9557,83285,157013,230741,...\}} = & 73728n + 9557 \\
\text{$d_{5,13} \in $ \{46421,193877,341333,488789,...\}} = & 147456n + 46421 \\
\text{$d_{5,14} \in $ \{267605,562517,857429,1152341,...\}} = & 294912n + 267605 \\
\text{$d_{5,15} \in $ \{120149,709973,1299797,1889621,...\}} = & 589824n + 120149 \\
\text{$d_{5,16} \in $ \{415061,1594709,2774357,3954005,...\}} = & 1179648n + 415061 \\
\text{$d_{5,17} \in $ \{1004885,3364181,5723477,8082773,...\}} = & 2359296n + 1004885 \\
\text{$d_{5,18} \in $ \{2184533,6903125,11621717,16340309,...\}} = & 4718592n + 2184533 \\
\vdots
\end{cases}
\rightarrow \frac{3d_5+1}{2^m} = odd
\end{equation}
}

When $i =6, r_6 = S(6) = 15$ (note that the index of V below is 15):
{\small
\begin{equation}
V_{15, m}(n) = 
\begin{cases}
\text{$V_{15,1} \in $ \{0,2,4,6,...\}} = & 2n \\
\text{$V_{15,2} \in $ \{1,5,9,13,...\}} = & 4n + 1 \\
\text{$V_{15,3} \in $ \{7,15,23,31,...\}} = & 8n + 7 \\
\text{$V_{15,4} \in $ \{3,19,35,51,...\}} = & 16n + 3 \\
\text{$V_{15,5} \in $ \{27,59,91,123,...\}} = & 32n + 27 \\
\text{$V_{15,6} \in $ \{43,107,171,235,...\}} = & 64n + 43 \\
\text{$V_{15,7} \in $ \{11,139,267,395,...\}} = & 128n + 11 \\
\text{$V_{15,8} \in $ \{203,459,715,971,...\}} = & 256n + 203 \\
\text{$V_{15,9} \in $ \{331,843,1355,1867,...\}} = & 512n + 331 \\
\text{$V_{15,10} \in $ \{587,1611,2635,3659,...\}} = & 1024n + 587 \\
\text{$V_{15,11} \in $ \{1099,3147,5195,7243,...\}} = & 2048n + 1099 \\
\text{$V_{15,12} \in $ \{75,4171,8267,12363,...\}} = & 4096n + 75 \\
\text{$V_{15,13} \in $ \{6219,14411,22603,30795,...\}} = & 8192n + 6219 \\
\text{$V_{15,14} \in $ \{2123,18507,34891,51275,...\}} = & 16384n + 2123 \\
\text{$V_{15,15} \in $ \{10315,43083,75851,108619,...\}} = & 32768n + 10315 \\
\text{$V_{15,16} \in $ \{59467,125003,190539,256075,...\}} = & 65536n + 59467 \\
\text{$V_{15,17} \in $ \{26699,157771,288843,419915,...\}} = & 131072n + 26699 \\
\text{$V_{15,18} \in $ \{92235,354379,616523,878667,...\}} = & 262144n + 92235 \\
\vdots
\end{cases}
\end{equation}
\begin{equation}
d_6 = 18V_{15, m}(n) +15 = 
\begin{cases}
\text{$d_{6,1} \in $ \{15,51,87,123,...\}} = & 36n + 15 \\
\text{$d_{6,2} \in $ \{33,105,177,249,...\}} = & 72n + 33 \\
\text{$d_{6,3} \in $ \{141,285,429,573,...\}} = & 144n + 141 \\
\text{$d_{6,4} \in $ \{69,357,645,933,...\}} = & 288n + 69 \\
\text{$d_{6,5} \in $ \{501,1077,1653,2229,...\}} = & 576n + 501 \\
\text{$d_{6,6} \in $ \{789,1941,3093,4245,...\}} = & 1152n + 789 \\
\text{$d_{6,7} \in $ \{213,2517,4821,7125,...\}} = & 2304n + 213 \\
\text{$d_{6,8} \in $ \{3669,8277,12885,17493,...\}} = & 4608n + 3669 \\
\text{$d_{6,9} \in $ \{5973,15189,24405,33621,...\}} = & 9216n + 5973 \\
\text{$d_{6,10} \in $ \{10581,29013,47445,65877,...\}} = & 18432n + 10581 \\
\text{$d_{6,11} \in $ \{19797,56661,93525,130389,...\}} = & 36864n + 19797 \\
\text{$d_{6,12} \in $ \{1365,75093,148821,222549,...\}} = & 73728n + 1365 \\
\text{$d_{6,13} \in $ \{111957,259413,406869,554325,...\}} = & 147456n + 111957 \\
\text{$d_{6,14} \in $ \{38229,333141,628053,922965,...\}} = & 294912n + 38229 \\
\text{$d_{6,15} \in $ \{185685,775509,1365333,1955157,...\}} = & 589824n + 185685 \\
\text{$d_{6,16} \in $ \{1070421,2250069,3429717,4609365,...\}} = & 1179648n + 1070421 \\
\text{$d_{6,17} \in $ \{480597,2839893,5199189,7558485,...\}} = & 2359296n + 480597 \\
\text{$d_{6,18} \in $ \{1660245,6378837,11097429,15816021,...\}} = & 4718592n + 1660245 \\
\vdots
\end{cases}
\rightarrow \frac{3d_6+1}{2^m} = odd
\end{equation}
}

When $i =7, r_7 = S(7) = 7$:
{\small
\begin{equation}
V_{7, m}(n) = 
\begin{cases}
\text{$V_{7,1} \in $ \{0,2,4,6,...\}} = & 2n \\
\text{$V_{7,2} \in $ \{1,5,9,13,...\}} = & 4n + 1 \\
\text{$V_{7,3} \in $ \{3,11,19,27,...\}} = & 8n + 3 \\
\text{$V_{7,4} \in $ \{7,23,39,55,...\}} = & 16n + 7 \\
\text{$V_{7,5} \in $ \{31,63,95,127,...\}} = & 32n + 31 \\
\text{$V_{7,6} \in $ \{15,79,143,207,...\}} = & 64n + 15 \\
\text{$V_{7,7} \in $ \{111,239,367,495,...\}} = & 128n + 111 \\
\text{$V_{7,8} \in $ \{175,431,687,943,...\}} = & 256n + 175 \\
\text{$V_{7,9} \in $ \{47,559,1071,1583,...\}} = & 512n + 47 \\
\text{$V_{7,10} \in $ \{815,1839,2863,3887,...\}} = & 1024n + 815 \\
\text{$V_{7,11} \in $ \{1327,3375,5423,7471,...\}} = & 2048n + 1327 \\
\text{$V_{7,12} \in $ \{2351,6447,10543,14639,...\}} = & 4096n + 2351 \\
\text{$V_{7,13} \in $ \{4399,12591,20783,28975,...\}} = & 8192n + 4399 \\
\text{$V_{7,14} \in $ \{303,16687,33071,49455,...\}} = & 16384n + 303 \\
\text{$V_{7,15} \in $ \{24879,57647,90415,123183,...\}} = & 32768n + 24879 \\
\text{$V_{7,16} \in $ \{8495,74031,139567,205103,...\}} = & 65536n + 8495 \\
\text{$V_{7,17} \in $ \{41263,172335,303407,434479,...\}} = & 131072n + 41263 \\
\text{$V_{7,18} \in $ \{237871,500015,762159,1024303,...\}} = & 262144n + 237871 \\
\vdots
\end{cases}
\end{equation}
\begin{equation}
d_7 = V_{7, m}(n) +7 =
\begin{cases}
\text{$d_{7,1} \in $ \{7,43,79,115,...\}} = & 36n + 7 \\
\text{$d_{7,2} \in $ \{25,97,169,241,...\}} = & 72n + 25 \\
\text{$d_{7,3} \in $ \{61,205,349,493,...\}} = & 144n + 61 \\
\text{$d_{7,4} \in $ \{133,421,709,997,...\}} = & 288n + 133 \\
\text{$d_{7,5} \in $ \{565,1141,1717,2293,...\}} = & 576n + 565 \\
\text{$d_{7,6} \in $ \{277,1429,2581,3733,...\}} = & 1152n + 277 \\
\text{$d_{7,7} \in $ \{2005,4309,6613,8917,...\}} = & 2304n + 2005 \\
\text{$d_{7,8} \in $ \{3157,7765,12373,16981,...\}} = & 4608n + 3157 \\
\text{$d_{7,9} \in $ \{853,10069,19285,28501,...\}} = & 9216n + 853 \\
\text{$d_{7,10} \in $ \{14677,33109,51541,69973,...\}} = & 18432n + 14677 \\
\text{$d_{7,11} \in $ \{23893,60757,97621,134485,...\}} = & 36864n + 23893 \\
\text{$d_{7,12} \in $ \{42325,116053,189781,263509,...\}} = & 73728n + 42325 \\
\text{$d_{7,13} \in $ \{79189,226645,374101,521557,...\}} = & 147456n + 79189 \\
\text{$d_{7,14} \in $ \{5461,300373,595285,890197,...\}} = & 294912n + 5461 \\
\text{$d_{7,15} \in $ \{447829,1037653,1627477,2217301,...\}} = & 589824n + 447829 \\
\text{$d_{7,16} \in $ \{152917,1332565,2512213,3691861,...\}} = & 1179648n + 152917 \\
\text{$d_{7,17} \in $ \{742741,3102037,5461333,7820629,...\}} = & 2359296n + 742741 \\
\text{$d_{7,18} \in $ \{4281685,9000277,13718869,18437461,...\}} = & 4718592n + 4281685 \\
\vdots
\end{cases}
\rightarrow \frac{3d_7+1}{2^m} = odd
\end{equation}
}

When $i =8, r_8 = S(8) = 11$ (note that the index of V below is 11):
{\small
\begin{equation}
V_{11, m}(n) = 
\begin{cases}
\text{$V_{11,1} \in $ \{0,2,4,6,...\}} = & 2n \\
\text{$V_{11,2} \in $ \{3,7,11,15,...\}} = & 4n + 3 \\
\text{$V_{11,3} \in $ \{1,9,17,25,...\}} = & 8n + 1 \\
\text{$V_{11,4} \in $ \{5,21,37,53,...\}} = & 16n + 5 \\
\text{$V_{11,5} \in $ \{13,45,77,109,...\}} = & 32n + 13 \\
\text{$V_{11,6} \in $ \{29,93,157,221,...\}} = & 64n + 29 \\
\text{$V_{11,7} \in $ \{125,253,381,509,...\}} = & 128n + 125 \\
\text{$V_{11,8} \in $ \{61,317,573,829,...\}} = & 256n + 61 \\
\text{$V_{11,9} \in $ \{445,957,1469,1981,...\}} = & 512n + 445 \\
\text{$V_{11,10} \in $ \{701,1725,2749,3773,...\}} = & 1024n + 701 \\
\text{$V_{11,11} \in $ \{189,2237,4285,6333,...\}} = & 2048n + 189 \\
\text{$V_{11,12} \in $ \{3261,7357,11453,15549,...\}} = & 4096n + 3261 \\
\text{$V_{11,13} \in $ \{5309,13501,21693,29885,...\}} = & 8192n + 5309 \\
\text{$V_{11,14} \in $ \{9405,25789,42173,58557,...\}} = & 16384n + 9405 \\
\text{$V_{11,15} \in $ \{17597,50365,83133,115901,...\}} = & 32768n + 17597 \\
\text{$V_{11,16} \in $ \{1213,66749,132285,197821,...\}} = & 65536n + 1213 \\
\text{$V_{11,17} \in $ \{99517,230589,361661,492733,...\}} = & 131072n + 99517 \\
\text{$V_{11,18} \in $ \{33981,296125,558269,820413,...\}} = & 262144n + 33981 \\
\vdots
\end{cases}
\end{equation}
\begin{equation}
d_8 = 18 V_{11, m}(n) +11
\begin{cases}
\text{$d_{8,1} \in $ \{11,47,83,119,...\}} = & 36n + 11 \\
\text{$d_{8,2} \in $ \{65,137,209,281,...\}} = & 72n + 65 \\
\text{$d_{8,3} \in $ \{29,173,317,461,...\}} = & 144n + 29 \\
\text{$d_{8,4} \in $ \{101,389,677,965,...\}} = & 288n + 101 \\
\text{$d_{8,5} \in $ \{245,821,1397,1973,...\}} = & 576n + 245 \\
\text{$d_{8,6} \in $ \{533,1685,2837,3989,...\}} = & 1152n + 533 \\
\text{$d_{8,7} \in $ \{2261,4565,6869,9173,...\}} = & 2304n + 2261 \\
\text{$d_{8,8} \in $ \{1109,5717,10325,14933,...\}} = & 4608n + 1109 \\
\text{$d_{8,9} \in $ \{8021,17237,26453,35669,...\}} = & 9216n + 8021 \\
\text{$d_{8,10} \in $ \{12629,31061,49493,67925,...\}} = & 18432n + 12629 \\
\text{$d_{8,11} \in $ \{3413,40277,77141,114005,...\}} = & 36864n + 3413 \\
\text{$d_{8,12} \in $ \{58709,132437,206165,279893,...\}} = & 73728n + 58709 \\
\text{$d_{8,13} \in $ \{95573,243029,390485,537941,...\}} = & 147456n + 95573 \\
\text{$d_{8,14} \in $ \{169301,464213,759125,1054037,...\}} = & 294912n + 169301 \\
\text{$d_{8,15} \in $ \{316757,906581,1496405,2086229,...\}} = & 589824n + 316757 \\
\text{$d_{8,16} \in $ \{21845,1201493,2381141,3560789,...\}} = & 1179648n + 21845 \\
\text{$d_{8,17} \in $ \{1791317,4150613,6509909,8869205,...\}} = & 2359296n + 1791317 \\
\text{$d_{8,18} \in $ \{611669,5330261,10048853,14767445,...\}} = & 4718592n + 611669 \\
\vdots
\end{cases}
\rightarrow \frac{3d_8+1}{2^m} = odd
\end{equation}
}

When $i =9, r_9 = S(9) = 9$:
{\small
\begin{equation}
V_{9, m}(n) = 
\begin{cases}
\text{$V_{9,1} \in $ \{1,3,5,7,...\}} = & 2n + 1 \\
\text{$V_{9,2} \in $ \{0,4,8,12,...\}} = & 4n \\
\text{$V_{9,3} \in $ \{2,10,18,26,...\}} = & 8n + 2 \\
\text{$V_{9,4} \in $ \{14,30,46,62,...\}} = & 16n + 14 \\
\text{$V_{9,5} \in $ \{6,38,70,102,...\}} = & 32n + 6 \\
\text{$V_{9,6} \in $ \{22,86,150,214,...\}} = & 64n + 22 \\
\text{$V_{9,7} \in $ \{54,182,310,438,...\}} = & 128n + 54 \\
\text{$V_{9,8} \in $ \{118,374,630,886,...\}} = & 256n + 118 \\
\text{$V_{9,9} \in $ \{502,1014,1526,2038,...\}} = & 512n + 502 \\
\text{$V_{9,10} \in $ \{246,1270,2294,3318,...\}} = & 1024n + 246 \\
\text{$V_{9,11} \in $ \{1782,3830,5878,7926,...\}} = & 2048n + 1782 \\
\text{$V_{9,12} \in $ \{2806,6902,10998,15094,...\}} = & 4096n + 2806 \\
\text{$V_{9,13} \in $ \{758,8950,17142,25334,...\}} = & 8192n + 758 \\
\text{$V_{9,14} \in $ \{13046,29430,45814,62198,...\}} = & 16384n + 13046 \\
\text{$V_{9,15} \in $ \{21238,54006,86774,119542,...\}} = & 32768n + 21238 \\
\text{$V_{9,16} \in $ \{37622,103158,168694,234230,...\}} = & 65536n + 37622 \\
\text{$V_{9,17} \in $ \{70390,201462,332534,463606,...\}} = & 131072n + 70390 \\
\text{$V_{9,18} \in $ \{4854,266998,529142,791286,...\}} = & 262144n + 4854 \\
\vdots
\end{cases}
\end{equation}
\begin{equation}
d_9 = 18V_{9, m}(n)+9 =
\begin{cases}
\text{$d_{9,1} \in $ \{27,63,99,135,...\}} = & 36n + 27 \\
\text{$d_{9,2} \in $ \{9,81,153,225,...\}} = & 72n + 9 \\
\text{$d_{9,3} \in $ \{45,189,333,477,...\}} = & 144n + 45 \\
\text{$d_{9,4} \in $ \{261,549,837,1125,...\}} = & 288n + 261 \\
\text{$d_{9,5} \in $ \{117,693,1269,1845,...\}} = & 576n + 117 \\
\text{$d_{9,6} \in $ \{405,1557,2709,3861,...\}} = & 1152n + 405 \\
\text{$d_{9,7} \in $ \{981,3285,5589,7893,...\}} = & 2304n + 981 \\
\text{$d_{9,8} \in $ \{2133,6741,11349,15957,...\}} = & 4608n + 2133 \\
\text{$d_{9,9} \in $ \{9045,18261,27477,36693,...\}} = & 9216n + 9045 \\
\text{$d_{9,10} \in $ \{4437,22869,41301,59733,...\}} = & 18432n + 4437 \\
\text{$d_{9,11} \in $ \{32085,68949,105813,142677,...\}} = & 36864n + 32085 \\
\text{$d_{9,12} \in $ \{50517,124245,197973,271701,...\}} = & 73728n + 50517 \\
\text{$d_{9,13} \in $ \{13653,161109,308565,456021,...\}} = & 147456n + 13653 \\
\text{$d_{9,14} \in $ \{234837,529749,824661,1119573,...\}} = & 294912n + 234837 \\
\text{$d_{9,15} \in $ \{382293,972117,1561941,2151765,...\}} = & 589824n + 382293 \\
\text{$d_{9,16} \in $ \{677205,1856853,3036501,4216149,...\}} = & 1179648n + 677205 \\
\text{$d_{9,17} \in $ \{1267029,3626325,5985621,8344917,...\}} = & 2359296n + 1267029 \\
\text{$d_{9,18} \in $ \{87381,4805973,9524565,14243157,...\}} = & 4718592n + 87381 \\
\vdots
\end{cases}
\rightarrow \frac{3d_9+1}{2^m} = odd
\end{equation}
}

The proofs of these $d_i$ statement is not difficult, for example 
\[
\begin{split}
\frac{3(36n+27)+1}{2^1}&=54n+41 \\
 \frac{3(72n+9)+1}{2^2}&=54n+7 \\ 
&\vdots \\
\frac{3(4718592n+87381)+1}{2^{18}}&=54n+1.
 \end{split}
 \]
\end{proof}

\begin{conjecture} \label{conj1} On the boundedness of $\frac{3d_i+1}{2^m}$.
One of the implications of theorem \ref{mainthm} is that the next odd number $d_{i_{next}}=\frac{3d_{i,m}+1}{2^m}$ is bounded between $54n$ and $54n+54$, i.e., $54n < d_{i_{next}} < 54(n+1)$  where $n$ is a variable dependent on $d_i$ (as also prescribed in the proposed Collatz map in \ref{CollMachinery1}). 
\end{conjecture} 

The implication of conj \ref{conj1} is $d_{i_{next}}$ must be congruent to residue $a$ modulo $54n$ where $ 1 \le a \le 53$ and  $a$ is odd. 

\begin{conjecture} \label{conj2}
These $d_i$ are values are the fundamental (principal) sets for all odd integers because they solely represent the sets from which all every odd number could be derived. 
\end{conjecture}

The proof of this last conjecture is beyond the scope of the objectives of this foundational paper.

As a result of the proposed conjecture in \ref{conj2} the map presented in \ref{CollMachinery1} is constructed from the $d_i$ formulations and proposed as the generalised Collatz based number system. This map \ref{CollMachinery1} consists of 3 major compartments: the top; middle; and the bottom sections representing the $d_i$, $3d_i+1$ and optimised $\frac{3d_i+1}{2^m}$ results, respectively.

Liewise, the map \ref{CollHeightMap1} consists of 3 major compartments: the top; middle; and the bottom sections representing the corresponding total stopping time functions of $d_i$, $3d_i+1$ and optimsed $\frac{3d_i+1}{2^m}$, respectively.


\begin{landscape}
{\fontsize{2.5}{4.0}\selectfont
\begin{equation}\label{CollMachinery1}
\begin{cases}
\textcolor{black}{{\bf Odd~d_1\downarrow}} \\
36n + 19*\\
\textcolor{blue}{ 72n + 1} \\
\textcolor{red}{ 144n + 109} \\
288n + 37\\
576n + 181\\
1152n + 1045\\
2304n + 469\\
\textcolor{blue}{ 4608n + 1621} \\
9216n + 3925\\
\textcolor{blue}{ 18432n + 8533} \\
36864n + 36181\\
73728n + 17749\\
\textcolor{blue}{ 147456n + 128341} \\
294912n + 202069\\
\textcolor{blue}{ 589824n + 54613} \\
1179648n + 939349\\
2359296n + 1529173\\
4718592n + 2708821\\
\vdots \\
\textcolor{black}{{\bf Even_1 \downarrow}} \\
108n + 58*\\
\textcolor{blue}{ 216n + 4} \\
\textcolor{red}{ 432n + 328} \\
864n + 112\\
1728n + 544\\
3456n + 3136\\
6912n + 1408\\
\textcolor{blue}{ 13824n + 4864} \\
27648n + 11776\\
\textcolor{blue}{ 55296n + 25600} \\
110592n + 108544\\
221184n + 53248\\
\textcolor{blue}{ 442368n + 385024} \\
884736n + 606208\\
\textcolor{blue}{ 1769472n + 163840} \\
3538944n + 2818048\\
7077888n + 4587520\\
14155776n + 8126464\\
\vdots \\
\textcolor{black}{{\bf Odd~d_{1_{next}} \downarrow}} \\
54n + 29*\\
\textcolor{blue}{ 54n + 1} \\
\textcolor{red}{ 54n + 41} \\
54n + 7\\
54n + 17\\
54n + 49\\
54n + 11\\
\textcolor{blue}{ 54n + 19} \\
54n + 23\\
\textcolor{blue}{ 54n + 25} \\
54n + 53\\
54n + 13\\
\textcolor{blue}{ 54n + 47} \\
54n + 37\\
\textcolor{blue}{ 54n + 5} \\
54n + 43\\
54n + 35\\
54n + 31\\
\vdots \\
\end{cases}
\begin{cases}
\textcolor{black}{{\bf Odd~d_2\downarrow}} \\
36n + 23*\\
72n + 41\\
144n + 77\\
\textcolor{blue}{ 288n + 5} \\
\textcolor{red}{ 576n + 437} \\
1152n + 149\\
2304n + 725\\
4608n + 4181\\
9216n + 1877\\
\textcolor{blue}{ 18432n + 6485} \\
36864n + 15701\\
\textcolor{blue}{ 73728n + 34133} \\
147456n + 144725\\
294912n + 70997\\
\textcolor{blue}{ 589824n + 513365} \\
1179648n + 808277\\
\textcolor{blue}{ 2359296n + 218453} \\
4718592n + 3757397\\
\vdots \\
\textcolor{black}{{\bf Even~_2 \downarrow}} \\
108n + 70*\\
216n + 124\\
432n + 232\\
\textcolor{blue}{ 864n + 16} \\
\textcolor{red}{ 1728n + 1312} \\
3456n + 448\\
6912n + 2176\\
13824n + 12544\\
27648n + 5632\\
\textcolor{blue}{ 55296n + 19456} \\
110592n + 47104\\
\textcolor{blue}{ 221184n + 102400} \\
442368n + 434176\\
884736n + 212992\\
\textcolor{blue}{ 1769472n + 1540096} \\
3538944n + 2424832\\
\textcolor{blue}{ 7077888n + 655360} \\
14155776n + 11272192\\
\vdots \\
\textcolor{black}{{\bf Odd~d_{2_{next}} \downarrow}} \\
54n + 35*\\
54n + 31\\
54n + 29\\
\textcolor{blue}{ 54n + 1} \\
\textcolor{red}{ 54n + 41} \\
54n + 7\\
54n + 17\\
54n + 49\\
54n + 11\\
\textcolor{blue}{ 54n + 19} \\
54n + 23\\
\textcolor{blue}{ 54n + 25} \\
54n + 53\\
54n + 13\\
\textcolor{blue}{ 54n + 47} \\
54n + 37\\
\textcolor{blue}{ 54n + 5} \\
54n + 43\\
\vdots \\
\end{cases}
\begin{cases}
\textcolor{black}{{\bf Odd~d_3\downarrow}} \\
\textcolor{blue}{ 36n + 3}* \\
72n + 57\\
144n + 93\\
288n + 165\\
576n + 309\\
\textcolor{blue}{ 1152n + 21} \\
\textcolor{red}{ 2304n + 1749} \\
4608n + 597\\
9216n + 2901\\
18432n + 16725\\
36864n + 7509\\
\textcolor{blue}{ 73728n + 25941} \\
147456n + 62805\\
\textcolor{blue}{ 294912n + 136533} \\
589824n + 578901\\
1179648n + 283989\\
\textcolor{blue}{ 2359296n + 2053461} \\
4718592n + 3233109\\
\vdots \\
\textcolor{black}{{\bf Even~_3 \downarrow}} \\
\textcolor{blue}{ 108n + 10}* \\
216n + 172\\
432n + 280\\
864n + 496\\
1728n + 928\\
\textcolor{blue}{ 3456n + 64} \\
\textcolor{red}{ 6912n + 5248} \\
13824n + 1792\\
27648n + 8704\\
55296n + 50176\\
110592n + 22528\\
\textcolor{blue}{ 221184n + 77824} \\
442368n + 188416\\
\textcolor{blue}{ 884736n + 409600} \\
1769472n + 1736704\\
3538944n + 851968\\
\textcolor{blue}{ 7077888n + 6160384} \\
14155776n + 9699328\\
\vdots \\
\textcolor{black}{{\bf Odd~d_{3_{next}} \downarrow}} \\
\textcolor{blue}{ 54n + 5}* \\
54n + 43\\
54n + 35\\
54n + 31\\
54n + 29\\
\textcolor{blue}{ 54n + 1} \\
\textcolor{red}{ 54n + 41} \\
54n + 7\\
54n + 17\\
54n + 49\\
54n + 11\\
\textcolor{blue}{ 54n + 19} \\
54n + 23\\
\textcolor{blue}{ 54n + 25} \\
54n + 53\\
54n + 13\\
\textcolor{blue}{ 54n + 47} \\
54n + 37\\
\vdots \\
\end{cases}
\begin{cases}
\textcolor{black}{{\bf Odd~d_4\downarrow}} \\
\textcolor{blue}{ 36n + 31}* \\
72n + 49\\
\textcolor{blue}{ 144n + 13} \\
288n + 229\\
576n + 373\\
1152n + 661\\
2304n + 1237\\
\textcolor{blue}{ 4608n + 85} \\
\textcolor{red}{ 9216n + 6997} \\
18432n + 2389\\
36864n + 11605\\
73728n + 66901\\
147456n + 30037\\
\textcolor{blue}{ 294912n + 103765} \\
589824n + 251221\\
\textcolor{blue}{ 1179648n + 546133} \\
2359296n + 2315605\\
4718592n + 1135957\\
\vdots \\
\textcolor{black}{{\bf Even~_4 \downarrow}} \\
\textcolor{blue}{ 108n + 94}* \\
216n + 148\\
\textcolor{blue}{ 432n + 40} \\
864n + 688\\
1728n + 1120\\
3456n + 1984\\
6912n + 3712\\
\textcolor{blue}{ 13824n + 256} \\
\textcolor{red}{ 27648n + 20992} \\
55296n + 7168\\
110592n + 34816\\
221184n + 200704\\
442368n + 90112\\
\textcolor{blue}{ 884736n + 311296} \\
1769472n + 753664\\
\textcolor{blue}{ 3538944n + 1638400} \\
7077888n + 6946816\\
14155776n + 3407872\\
\vdots \\
\textcolor{black}{{\bf Odd~d_{4_{next}} \downarrow}} \\
\textcolor{blue}{ 54n + 47}* \\
54n + 37\\
\textcolor{blue}{ 54n + 5} \\
54n + 43\\
54n + 35\\
54n + 31\\
54n + 29\\
\textcolor{blue}{ 54n + 1} \\
\textcolor{red}{ 54n + 41} \\
54n + 7\\
54n + 17\\
54n + 49\\
54n + 11\\
\textcolor{blue}{ 54n + 19} \\
54n + 23\\
\textcolor{blue}{ 54n + 25} \\
54n + 53\\
54n + 13\\
\vdots \\
\end{cases}
\begin{cases}
\textcolor{black}{{\bf Odd~d_5\downarrow}} \\
36n + 35*\\
72n + 17\\
\textcolor{blue}{ 144n + 125} \\
288n + 197\\
\textcolor{blue}{ 576n + 53} \\
1152n + 917\\
2304n + 1493\\
4608n + 2645\\
9216n + 4949\\
\textcolor{blue}{ 18432n + 341} \\
\textcolor{red}{ 36864n + 27989} \\
73728n + 9557\\
147456n + 46421\\
294912n + 267605\\
589824n + 120149\\
\textcolor{blue}{ 1179648n + 415061} \\
2359296n + 1004885\\
4718592n + 2184533\\
\vdots \\
\textcolor{black}{{\bf Even~_5 \downarrow}} \\
108n + 106*\\
216n + 52\\
\textcolor{blue}{ 432n + 376} \\
864n + 592\\
\textcolor{blue}{ 1728n + 160} \\
3456n + 2752\\
6912n + 4480\\
13824n + 7936\\
27648n + 14848\\
\textcolor{blue}{ 55296n + 1024} \\
\textcolor{red}{ 110592n + 83968} \\
221184n + 28672\\
442368n + 139264\\
884736n + 802816\\
1769472n + 360448\\
\textcolor{blue}{ 3538944n + 1245184} \\
7077888n + 3014656\\
14155776n + 6553600\\
\vdots \\
\textcolor{black}{{\bf Odd~d_{5_{next}} \downarrow}} \\
54n + 53*\\
54n + 13\\
\textcolor{blue}{ 54n + 47} \\
54n + 37\\
\textcolor{blue}{ 54n + 5} \\
54n + 43\\
54n + 35\\
54n + 31\\
54n + 29\\
\textcolor{blue}{ 54n + 1} \\
\textcolor{red}{ 54n + 41} \\
54n + 7\\
54n + 17\\
54n + 49\\
54n + 11\\
\textcolor{blue}{ 54n + 19} \\
54n + 23\\
54n + 25\\
\vdots \\
\end{cases}
\begin{cases}
\textcolor{black}{{\bf Odd~d_6\downarrow}} \\
36n + 15*\\
\textcolor{blue}{ 72n + 33} \\
144n + 141\\
288n + 69\\
\textcolor{blue}{ 576n + 501} \\
1152n + 789\\
\textcolor{blue}{ 2304n + 213} \\
4608n + 3669\\
9216n + 5973\\
18432n + 10581\\
36864n + 19797\\
\textcolor{blue}{ 73728n + 1365} \\
\textcolor{red}{ 147456n + 111957} \\
294912n + 38229\\
589824n + 185685\\
1179648n + 1070421\\
2359296n + 480597\\
4718592n + 1660245\\
\vdots \\
\textcolor{black}{{\bf Even~_6 \downarrow}} \\
108n + 46*\\
\textcolor{blue}{ 216n + 100} \\
432n + 424\\
864n + 208\\
\textcolor{blue}{ 1728n + 1504} \\
3456n + 2368\\
\textcolor{blue}{ 6912n + 640} \\
13824n + 11008\\
27648n + 17920\\
55296n + 31744\\
110592n + 59392\\
\textcolor{blue}{ 221184n + 4096} \\
\textcolor{red}{ 442368n + 335872} \\
884736n + 114688\\
1769472n + 557056\\
3538944n + 3211264\\
7077888n + 1441792\\
14155776n + 4980736\\
\vdots \\
\textcolor{black}{{\bf Odd~d_{6_{next}} \downarrow}} \\
54n + 23*\\
\textcolor{blue}{ 54n + 25} \\
54n + 53\\
54n + 13\\
\textcolor{blue}{ 54n + 47} \\
54n + 37\\
\textcolor{blue}{ 54n + 5} \\
54n + 43\\
54n + 35\\
54n + 31\\
54n + 29\\
\textcolor{blue}{ 54n + 1} \\
\textcolor{red}{ 54n + 41} \\
54n + 7\\
54n + 17\\
54n + 49\\
54n + 11\\
54n + 19\\
\vdots \\
\end{cases}
\begin{cases}
\textcolor{black}{{\bf Odd~d_7\downarrow}} \\
36n + 7*\\
\textcolor{blue}{ 72n + 25} \\
144n + 61\\
\textcolor{blue}{ 288n + 133} \\
576n + 565\\
1152n + 277\\
\textcolor{blue}{ 2304n + 2005} \\
4608n + 3157\\
\textcolor{blue}{ 9216n + 853} \\
18432n + 14677\\
36864n + 23893\\
73728n + 42325\\
147456n + 79189\\
\textcolor{blue}{ 294912n + 5461} \\
\textcolor{red}{ 589824n + 447829} \\
1179648n + 152917\\
2359296n + 742741\\
4718592n + 4281685\\
\vdots \\
\textcolor{black}{{\bf Even~_7 \downarrow}} \\
108n + 22*\\
\textcolor{blue}{ 216n + 76} \\
432n + 184\\
\textcolor{blue}{ 864n + 400} \\
1728n + 1696\\
3456n + 832\\
\textcolor{blue}{ 6912n + 6016} \\
13824n + 9472\\
\textcolor{blue}{ 27648n + 2560} \\
55296n + 44032\\
110592n + 71680\\
221184n + 126976\\
442368n + 237568\\
\textcolor{blue}{ 884736n + 16384} \\
\textcolor{red}{ 1769472n + 1343488} \\
3538944n + 458752\\
7077888n + 2228224\\
14155776n + 12845056\\
\vdots \\
\textcolor{black}{{\bf Odd~d_{7_{next}} \downarrow}} \\
54n + 11*\\
\textcolor{blue}{ 54n + 19} \\
54n + 23\\
\textcolor{blue}{ 54n + 25} \\
54n + 53\\
54n + 13\\
\textcolor{blue}{ 54n + 47} \\
54n + 37\\
\textcolor{blue}{ 54n + 5} \\
54n + 43\\
54n + 35\\
54n + 31\\
54n + 29\\
\textcolor{blue}{ 54n + 1} \\
\textcolor{red}{ 54n + 41} \\
54n + 7\\
54n + 17\\
54n + 49\\
\vdots \\
\end{cases}
\begin{cases}
\textcolor{black}{{\bf Odd~d_8\downarrow}} \\
36n + 11*\\
72n + 65\\
144n + 29\\
\textcolor{blue}{ 288n + 101} \\
576n + 245\\
\textcolor{blue}{ 1152n + 533} \\
2304n + 2261\\
4608n + 1109\\
\textcolor{blue}{ 9216n + 8021} \\
18432n + 12629\\
\textcolor{blue}{ 36864n + 3413} \\
73728n + 58709\\
147456n + 95573\\
294912n + 169301\\
589824n + 316757\\
\textcolor{blue}{ 1179648n + 21845} \\
\textcolor{red}{ 2359296n + 1791317} \\
4718592n + 611669\\
\vdots \\
\textcolor{black}{{\bf Even~_8 \downarrow}} \\
108n + 34*\\
216n + 196\\
432n + 88\\
\textcolor{blue}{ 864n + 304} \\
1728n + 736\\
\textcolor{blue}{ 3456n + 1600} \\
6912n + 6784\\
13824n + 3328\\
\textcolor{blue}{ 27648n + 24064} \\
55296n + 37888\\
\textcolor{blue}{ 110592n + 10240} \\
221184n + 176128\\
442368n + 286720\\
884736n + 507904\\
1769472n + 950272\\
\textcolor{blue}{ 3538944n + 65536} \\
\textcolor{red}{ 7077888n + 5373952} \\
14155776n + 1835008\\
\vdots \\
\textcolor{black}{{\bf Odd~d_{8_{next}} \downarrow}} \\
54n + 17*\\
54n + 49\\
54n + 11\\
\textcolor{blue}{ 54n + 19} \\
54n + 23\\
\textcolor{blue}{ 54n + 25} \\
54n + 53\\
54n + 13\\
\textcolor{blue}{ 54n + 47} \\
54n + 37\\
\textcolor{blue}{ 54n + 5} \\
54n + 43\\
54n + 35\\
54n + 31\\
54n + 29\\
\textcolor{blue}{ 54n + 1} \\
\textcolor{red}{ 54n + 41} \\
54n + 7\\
\vdots \\
\end{cases}
\begin{cases}
\textcolor{black}{{\bf Odd~d_9\downarrow}} \\
\textcolor{red}{ 36n + 27}* \\
72n + 9\\
144n + 45\\
288n + 261\\
576n + 117\\
\textcolor{blue}{ 1152n + 405} \\
2304n + 981\\
\textcolor{blue}{ 4608n + 2133} \\
9216n + 9045\\
18432n + 4437\\
\textcolor{blue}{ 36864n + 32085} \\
73728n + 50517\\
\textcolor{blue}{ 147456n + 13653} \\
294912n + 234837\\
589824n + 382293\\
1179648n + 677205\\
2359296n + 1267029\\
\textcolor{blue}{ 4718592n + 87381} \\
\vdots \\
\textcolor{black}{{\bf Even~_9 \downarrow}} \\
\textcolor{red}{ 108n + 82}* \\
216n + 28\\
432n + 136\\
864n + 784\\
1728n + 352\\
\textcolor{blue}{ 3456n + 1216} \\
6912n + 2944\\
\textcolor{blue}{ 13824n + 6400} \\
27648n + 27136\\
55296n + 13312\\
\textcolor{blue}{ 110592n + 96256} \\
221184n + 151552\\
\textcolor{blue}{ 442368n + 40960} \\
884736n + 704512\\
1769472n + 1146880\\
3538944n + 2031616\\
7077888n + 3801088\\
\textcolor{blue}{ 14155776n + 262144} \\
\vdots \\
\textcolor{black}{{\bf Odd~d_{9_{next}} \downarrow}} \\
\textcolor{red}{ 54n + 41}* \\
54n + 7\\
54n + 17\\
54n + 49\\
54n + 11\\
\textcolor{blue}{ 54n + 19} \\
54n + 23\\
\textcolor{blue}{ 54n + 25} \\
54n + 53\\
54n + 13\\
\textcolor{blue}{ 54n + 47} \\
54n + 37\\
\textcolor{blue}{ 54n + 5} \\
54n + 43\\
54n + 35\\
54n + 31\\
54n + 29\\
\textcolor{blue}{ 54n + 1} \\
\vdots \\
\end{cases}
\footnote{All rights reserved. $\copyright$ Michael A. Idowu, 2014.}
\end{equation}
}
\end{landscape}

\begin{landscape}
{\fontsize{2.5}{4.0}\selectfont
\begin{equation}\label{CollHeightMap1}
\begin{cases}
\textcolor{black}{{\bf Odd~d_1\downarrow}} \\
\sigma_{\infty}(54n+29)+2\\
\textcolor{blue}{ \sigma_{\infty}(54n+1)+3} \\
\textcolor{red}{ \sigma_{\infty}(54n+41)+4} \\
\sigma_{\infty}(54n+7)+5\\
\sigma_{\infty}(54n+17)+6\\
\sigma_{\infty}(54n+49)+7\\
\sigma_{\infty}(54n+11)+8\\
\textcolor{blue}{ \sigma_{\infty}(54n+19)+9} \\
\sigma_{\infty}(54n+23)+10\\
\textcolor{blue}{ \sigma_{\infty}(54n+25)+11} \\
\sigma_{\infty}(54n+53)+12\\
\sigma_{\infty}(54n+13)+13\\
\textcolor{blue}{ \sigma_{\infty}(54n+47)+14} \\
\sigma_{\infty}(54n+37)+15\\
\textcolor{blue}{ \sigma_{\infty}(54n+5)+16} \\
\sigma_{\infty}(54n+43)+17\\
\sigma_{\infty}(54n+35)+18\\
\sigma_{\infty}(54n+31)+19\\
\vdots \\
\textcolor{black}{{\bf Even_1 \downarrow}} \\
\sigma_{\infty}(54n+29)+1\\
\textcolor{blue}{ \sigma_{\infty}(54n+1)+2} \\
\textcolor{red}{ \sigma_{\infty}(54n+41)+3} \\
\sigma_{\infty}(54n+7)+4\\
\sigma_{\infty}(54n+17)+5\\
\sigma_{\infty}(54n+49)+6\\
\sigma_{\infty}(54n+11)+7\\
\textcolor{blue}{ \sigma_{\infty}(54n+19)+8} \\
\sigma_{\infty}(54n+23)+9\\
\textcolor{blue}{ \sigma_{\infty}(54n+25)+10} \\
\sigma_{\infty}(54n+53)+11\\
\sigma_{\infty}(54n+13)+12\\
\textcolor{blue}{ \sigma_{\infty}(54n+47)+13} \\
\sigma_{\infty}(54n+37)+14\\
\textcolor{blue}{ \sigma_{\infty}(54n+5)+15} \\
\sigma_{\infty}(54n+43)+16\\
\sigma_{\infty}(54n+35)+17\\
\sigma_{\infty}(54n+31)+18\\
\vdots \\
\textcolor{black}{{\bf Odd~d_{1_{next}} \downarrow}} \\
\sigma_{\infty}(54n+29) \\
\textcolor{blue}{ \sigma_{\infty}(54n+1) } \\
\textcolor{red}{ \sigma_{\infty}(54n+41) } \\
\sigma_{\infty}(54n+7) \\
\sigma_{\infty}(54n+17) \\
\sigma_{\infty}(54n+49) \\
\sigma_{\infty}(54n+11) \\
\textcolor{blue}{ \sigma_{\infty}(54n+19) } \\
\sigma_{\infty}(54n+23) \\
\textcolor{blue}{ \sigma_{\infty}(54n+25) } \\
\sigma_{\infty}(54n+53) \\
\sigma_{\infty}(54n+13) \\
\textcolor{blue}{ \sigma_{\infty}(54n+47) } \\
\sigma_{\infty}(54n+37) \\
\textcolor{blue}{ \sigma_{\infty}(54n+5) } \\
\sigma_{\infty}(54n+43) \\
\sigma_{\infty}(54n+35) \\
\sigma_{\infty}(54n+31) \\
\vdots \\
\end{cases}
\begin{cases}
\textcolor{black}{{\bf Odd~d_2\downarrow}} \\
\sigma_{\infty}(54n+35)+2\\
\sigma_{\infty}(54n+31)+3\\
\sigma_{\infty}(54n+29)+4\\
\textcolor{blue}{ \sigma_{\infty}(54n+1)+5} \\
\textcolor{red}{ \sigma_{\infty}(54n+41)+6} \\
\sigma_{\infty}(54n+7)+7\\
\sigma_{\infty}(54n+17)+8\\
\sigma_{\infty}(54n+49)+9\\
\sigma_{\infty}(54n+11)+10\\
\textcolor{blue}{ \sigma_{\infty}(54n+19)+11} \\
\sigma_{\infty}(54n+23)+12\\
\textcolor{blue}{ \sigma_{\infty}(54n+25)+13} \\
\sigma_{\infty}(54n+53)+14\\
\sigma_{\infty}(54n+13)+15\\
\textcolor{blue}{ \sigma_{\infty}(54n+47)+16} \\
\sigma_{\infty}(54n+37)+17\\
\textcolor{blue}{ \sigma_{\infty}(54n+5)+18} \\
\sigma_{\infty}(54n+43)+19\\
\vdots \\
\textcolor{black}{{\bf Even_2 \downarrow}} \\
\sigma_{\infty}(54n+35)+1\\
\sigma_{\infty}(54n+31)+2\\
\sigma_{\infty}(54n+29)+3\\
\textcolor{blue}{ \sigma_{\infty}(54n+1)+4} \\
\textcolor{red}{ \sigma_{\infty}(54n+41)+5} \\
\sigma_{\infty}(54n+7)+6\\
\sigma_{\infty}(54n+17)+7\\
\sigma_{\infty}(54n+49)+8\\
\sigma_{\infty}(54n+11)+9\\
\textcolor{blue}{ \sigma_{\infty}(54n+19)+10} \\
\sigma_{\infty}(54n+23)+11\\
\textcolor{blue}{ \sigma_{\infty}(54n+25)+12} \\
\sigma_{\infty}(54n+53)+13\\
\sigma_{\infty}(54n+13)+14\\
\textcolor{blue}{ \sigma_{\infty}(54n+47)+15} \\
\sigma_{\infty}(54n+37)+16\\
\textcolor{blue}{ \sigma_{\infty}(54n+5)+17} \\
\sigma_{\infty}(54n+43)+18\\
\vdots \\
\textcolor{black}{{\bf Odd~d_{2_{next}} \downarrow}} \\
\sigma_{\infty}(54n+35) \\
\sigma_{\infty}(54n+31) \\
\sigma_{\infty}(54n+29) \\
\textcolor{blue}{ \sigma_{\infty}(54n+1) } \\
\textcolor{red}{ \sigma_{\infty}(54n+41) } \\
\sigma_{\infty}(54n+7) \\
\sigma_{\infty}(54n+17) \\
\sigma_{\infty}(54n+49) \\
\sigma_{\infty}(54n+11) \\
\textcolor{blue}{ \sigma_{\infty}(54n+19) } \\
\sigma_{\infty}(54n+23) \\
\textcolor{blue}{ \sigma_{\infty}(54n+25) } \\
\sigma_{\infty}(54n+53) \\
\sigma_{\infty}(54n+13) \\
\textcolor{blue}{ \sigma_{\infty}(54n+47) } \\
\sigma_{\infty}(54n+37) \\
\textcolor{blue}{ \sigma_{\infty}(54n+5) } \\
\sigma_{\infty}(54n+43) \\
\vdots \\
\end{cases}
\begin{cases}
\textcolor{black}{{\bf Odd~d_3\downarrow}} \\
\textcolor{blue}{ \sigma_{\infty}(54n+5)+2} \\
\sigma_{\infty}(54n+43)+3\\
\sigma_{\infty}(54n+35)+4\\
\sigma_{\infty}(54n+31)+5\\
\sigma_{\infty}(54n+29)+6\\
\textcolor{blue}{ \sigma_{\infty}(54n+1)+7} \\
\textcolor{red}{ \sigma_{\infty}(54n+41)+8} \\
\sigma_{\infty}(54n+7)+9\\
\sigma_{\infty}(54n+17)+10\\
\sigma_{\infty}(54n+49)+11\\
\sigma_{\infty}(54n+11)+12\\
\textcolor{blue}{ \sigma_{\infty}(54n+19)+13} \\
\sigma_{\infty}(54n+23)+14\\
\textcolor{blue}{ \sigma_{\infty}(54n+25)+15} \\
\sigma_{\infty}(54n+53)+16\\
\sigma_{\infty}(54n+13)+17\\
\textcolor{blue}{ \sigma_{\infty}(54n+47)+18} \\
\sigma_{\infty}(54n+37)+19\\
\vdots \\
\textcolor{black}{{\bf Even_3 \downarrow}} \\
\textcolor{blue}{ \sigma_{\infty}(54n+5)+1} \\
\sigma_{\infty}(54n+43)+2\\
\sigma_{\infty}(54n+35)+3\\
\sigma_{\infty}(54n+31)+4\\
\sigma_{\infty}(54n+29)+5\\
\textcolor{blue}{ \sigma_{\infty}(54n+1)+6} \\
\textcolor{red}{ \sigma_{\infty}(54n+41)+7} \\
\sigma_{\infty}(54n+7)+8\\
\sigma_{\infty}(54n+17)+9\\
\sigma_{\infty}(54n+49)+10\\
\sigma_{\infty}(54n+11)+11\\
\textcolor{blue}{ \sigma_{\infty}(54n+19)+12} \\
\sigma_{\infty}(54n+23)+13\\
\textcolor{blue}{ \sigma_{\infty}(54n+25)+14} \\
\sigma_{\infty}(54n+53)+15\\
\sigma_{\infty}(54n+13)+16\\
\textcolor{blue}{ \sigma_{\infty}(54n+47)+17} \\
\sigma_{\infty}(54n+37)+18\\
\vdots \\
\textcolor{black}{{\bf Odd~d_{3_{next}} \downarrow}} \\
\textcolor{blue}{ \sigma_{\infty}(54n+5) } \\
\sigma_{\infty}(54n+43) \\
\sigma_{\infty}(54n+35) \\
\sigma_{\infty}(54n+31) \\
\sigma_{\infty}(54n+29) \\
\textcolor{blue}{ \sigma_{\infty}(54n+1) } \\
\textcolor{red}{ \sigma_{\infty}(54n+41) } \\
\sigma_{\infty}(54n+7) \\
\sigma_{\infty}(54n+17) \\
\sigma_{\infty}(54n+49) \\
\sigma_{\infty}(54n+11) \\
\textcolor{blue}{ \sigma_{\infty}(54n+19) } \\
\sigma_{\infty}(54n+23) \\
\textcolor{blue}{ \sigma_{\infty}(54n+25) } \\
\sigma_{\infty}(54n+53) \\
\sigma_{\infty}(54n+13) \\
\textcolor{blue}{ \sigma_{\infty}(54n+47) } \\
\sigma_{\infty}(54n+37) \\
\vdots \\
\end{cases}
\begin{cases}
\textcolor{black}{{\bf Odd~d_4\downarrow}} \\
\textcolor{blue}{ \sigma_{\infty}(54n+47)+2} \\
\sigma_{\infty}(54n+37)+3\\
\textcolor{blue}{ \sigma_{\infty}(54n+5)+4} \\
\sigma_{\infty}(54n+43)+5\\
\sigma_{\infty}(54n+35)+6\\
\sigma_{\infty}(54n+31)+7\\
\sigma_{\infty}(54n+29)+8\\
\textcolor{blue}{ \sigma_{\infty}(54n+1)+9} \\
\textcolor{red}{ \sigma_{\infty}(54n+41)+10} \\
\sigma_{\infty}(54n+7)+11\\
\sigma_{\infty}(54n+17)+12\\
\sigma_{\infty}(54n+49)+13\\
\sigma_{\infty}(54n+11)+14\\
\textcolor{blue}{ \sigma_{\infty}(54n+19)+15} \\
\sigma_{\infty}(54n+23)+16\\
\textcolor{blue}{ \sigma_{\infty}(54n+25)+17} \\
\sigma_{\infty}(54n+53)+18\\
\sigma_{\infty}(54n+13)+19\\
\vdots \\
\textcolor{black}{{\bf Even_4 \downarrow}} \\
\textcolor{blue}{ \sigma_{\infty}(54n+47)+1} \\
\sigma_{\infty}(54n+37)+2\\
\textcolor{blue}{ \sigma_{\infty}(54n+5)+3} \\
\sigma_{\infty}(54n+43)+4\\
\sigma_{\infty}(54n+35)+5\\
\sigma_{\infty}(54n+31)+6\\
\sigma_{\infty}(54n+29)+7\\
\textcolor{blue}{ \sigma_{\infty}(54n+1)+8} \\
\textcolor{red}{ \sigma_{\infty}(54n+41)+9} \\
\sigma_{\infty}(54n+7)+10\\
\sigma_{\infty}(54n+17)+11\\
\sigma_{\infty}(54n+49)+12\\
\sigma_{\infty}(54n+11)+13\\
\textcolor{blue}{ \sigma_{\infty}(54n+19)+14} \\
\sigma_{\infty}(54n+23)+15\\
\textcolor{blue}{ \sigma_{\infty}(54n+25)+16} \\
\sigma_{\infty}(54n+53)+17\\
\sigma_{\infty}(54n+13)+18\\
\vdots \\
\textcolor{black}{{\bf Odd~d_{4_{next}} \downarrow}} \\
\textcolor{blue}{ \sigma_{\infty}(54n+47) } \\
\sigma_{\infty}(54n+37) \\
\textcolor{blue}{ \sigma_{\infty}(54n+5) } \\
\sigma_{\infty}(54n+43) \\
\sigma_{\infty}(54n+35) \\
\sigma_{\infty}(54n+31) \\
\sigma_{\infty}(54n+29) \\
\textcolor{blue}{ \sigma_{\infty}(54n+1) } \\
\textcolor{red}{ \sigma_{\infty}(54n+41) } \\
\sigma_{\infty}(54n+7) \\
\sigma_{\infty}(54n+17) \\
\sigma_{\infty}(54n+49) \\
\sigma_{\infty}(54n+11) \\
\textcolor{blue}{ \sigma_{\infty}(54n+19) } \\
\sigma_{\infty}(54n+23) \\
\textcolor{blue}{ \sigma_{\infty}(54n+25) } \\
\sigma_{\infty}(54n+53) \\
\sigma_{\infty}(54n+13) \\
\vdots \\
\end{cases}
\begin{cases}
\textcolor{black}{{\bf Odd~d_5\downarrow}} \\
\sigma_{\infty}(54n+53)+2\\
\sigma_{\infty}(54n+13)+3\\
\textcolor{blue}{ \sigma_{\infty}(54n+47)+4} \\
\sigma_{\infty}(54n+37)+5\\
\textcolor{blue}{ \sigma_{\infty}(54n+5)+6} \\
\sigma_{\infty}(54n+43)+7\\
\sigma_{\infty}(54n+35)+8\\
\sigma_{\infty}(54n+31)+9\\
\sigma_{\infty}(54n+29)+10\\
\textcolor{blue}{ \sigma_{\infty}(54n+1)+11} \\
\textcolor{red}{ \sigma_{\infty}(54n+41)+12} \\
\sigma_{\infty}(54n+7)+13\\
\sigma_{\infty}(54n+17)+14\\
\sigma_{\infty}(54n+49)+15\\
\sigma_{\infty}(54n+11)+16\\
\textcolor{blue}{ \sigma_{\infty}(54n+19)+17} \\
\sigma_{\infty}(54n+23)+18\\
\sigma_{\infty}(54n+25)+19\\
\vdots \\
\textcolor{black}{{\bf Even_5 \downarrow}} \\
\sigma_{\infty}(54n+53)+1\\
\sigma_{\infty}(54n+13)+2\\
\textcolor{blue}{ \sigma_{\infty}(54n+47)+3} \\
\sigma_{\infty}(54n+37)+4\\
\textcolor{blue}{ \sigma_{\infty}(54n+5)+5} \\
\sigma_{\infty}(54n+43)+6\\
\sigma_{\infty}(54n+35)+7\\
\sigma_{\infty}(54n+31)+8\\
\sigma_{\infty}(54n+29)+9\\
\textcolor{blue}{ \sigma_{\infty}(54n+1)+10} \\
\textcolor{red}{ \sigma_{\infty}(54n+41)+11} \\
\sigma_{\infty}(54n+7)+12\\
\sigma_{\infty}(54n+17)+13\\
\sigma_{\infty}(54n+49)+14\\
\sigma_{\infty}(54n+11)+15\\
\textcolor{blue}{ \sigma_{\infty}(54n+19)+16} \\
\sigma_{\infty}(54n+23)+17\\
\sigma_{\infty}(54n+25)+18\\
\vdots \\
\textcolor{black}{{\bf Odd~d_{5_{next}} \downarrow}} \\
\sigma_{\infty}(54n+53) \\
\sigma_{\infty}(54n+13) \\
\textcolor{blue}{ \sigma_{\infty}(54n+47) } \\
\sigma_{\infty}(54n+37) \\
\textcolor{blue}{ \sigma_{\infty}(54n+5) } \\
\sigma_{\infty}(54n+43) \\
\sigma_{\infty}(54n+35) \\
\sigma_{\infty}(54n+31) \\
\sigma_{\infty}(54n+29) \\
\textcolor{blue}{ \sigma_{\infty}(54n+1) } \\
\textcolor{red}{ \sigma_{\infty}(54n+41) } \\
\sigma_{\infty}(54n+7) \\
\sigma_{\infty}(54n+17) \\
\sigma_{\infty}(54n+49) \\
\sigma_{\infty}(54n+11) \\
\textcolor{blue}{ \sigma_{\infty}(54n+19) } \\
\sigma_{\infty}(54n+23) \\
\sigma_{\infty}(54n+25) \\
\vdots \\
\end{cases}
\begin{cases}
\textcolor{black}{{\bf Odd~d_6\downarrow}} \\
\sigma_{\infty}(54n+23)+2\\
\textcolor{blue}{ \sigma_{\infty}(54n+25)+3} \\
\sigma_{\infty}(54n+53)+4\\
\sigma_{\infty}(54n+13)+5\\
\textcolor{blue}{ \sigma_{\infty}(54n+47)+6} \\
\sigma_{\infty}(54n+37)+7\\
\textcolor{blue}{ \sigma_{\infty}(54n+5)+8} \\
\sigma_{\infty}(54n+43)+9\\
\sigma_{\infty}(54n+35)+10\\
\sigma_{\infty}(54n+31)+11\\
\sigma_{\infty}(54n+29)+12\\
\textcolor{blue}{ \sigma_{\infty}(54n+1)+13} \\
\textcolor{red}{ \sigma_{\infty}(54n+41)+14} \\
\sigma_{\infty}(54n+7)+15\\
\sigma_{\infty}(54n+17)+16\\
\sigma_{\infty}(54n+49)+17\\
\sigma_{\infty}(54n+11)+18\\
\sigma_{\infty}(54n+19)+19\\
\vdots \\
\textcolor{black}{{\bf Even_6 \downarrow}} \\
\sigma_{\infty}(54n+23)+1\\
\textcolor{blue}{ \sigma_{\infty}(54n+25)+2} \\
\sigma_{\infty}(54n+53)+3\\
\sigma_{\infty}(54n+13)+4\\
\textcolor{blue}{ \sigma_{\infty}(54n+47)+5} \\
\sigma_{\infty}(54n+37)+6\\
\textcolor{blue}{ \sigma_{\infty}(54n+5)+7} \\
\sigma_{\infty}(54n+43)+8\\
\sigma_{\infty}(54n+35)+9\\
\sigma_{\infty}(54n+31)+10\\
\sigma_{\infty}(54n+29)+11\\
\textcolor{blue}{ \sigma_{\infty}(54n+1)+12} \\
\textcolor{red}{ \sigma_{\infty}(54n+41)+13} \\
\sigma_{\infty}(54n+7)+14\\
\sigma_{\infty}(54n+17)+15\\
\sigma_{\infty}(54n+49)+16\\
\sigma_{\infty}(54n+11)+17\\
\sigma_{\infty}(54n+19)+18\\
\vdots \\
\textcolor{black}{{\bf Odd~d_{6_{next}} \downarrow}} \\
\sigma_{\infty}(54n+23) \\
\textcolor{blue}{ \sigma_{\infty}(54n+25) } \\
\sigma_{\infty}(54n+53) \\
\sigma_{\infty}(54n+13) \\
\textcolor{blue}{ \sigma_{\infty}(54n+47) } \\
\sigma_{\infty}(54n+37) \\
\textcolor{blue}{ \sigma_{\infty}(54n+5) } \\
\sigma_{\infty}(54n+43) \\
\sigma_{\infty}(54n+35) \\
\sigma_{\infty}(54n+31) \\
\sigma_{\infty}(54n+29) \\
\textcolor{blue}{ \sigma_{\infty}(54n+1) } \\
\textcolor{red}{ \sigma_{\infty}(54n+41) } \\
\sigma_{\infty}(54n+7) \\
\sigma_{\infty}(54n+17) \\
\sigma_{\infty}(54n+49) \\
\sigma_{\infty}(54n+11) \\
\sigma_{\infty}(54n+19) \\
\vdots \\
\end{cases}
\begin{cases}
\textcolor{black}{{\bf Odd~d_7\downarrow}} \\
\sigma_{\infty}(54n+11)+2\\
\textcolor{blue}{ \sigma_{\infty}(54n+19)+3} \\
\sigma_{\infty}(54n+23)+4\\
\textcolor{blue}{ \sigma_{\infty}(54n+25)+5} \\
\sigma_{\infty}(54n+53)+6\\
\sigma_{\infty}(54n+13)+7\\
\textcolor{blue}{ \sigma_{\infty}(54n+47)+8} \\
\sigma_{\infty}(54n+37)+9\\
\textcolor{blue}{ \sigma_{\infty}(54n+5)+10} \\
\sigma_{\infty}(54n+43)+11\\
\sigma_{\infty}(54n+35)+12\\
\sigma_{\infty}(54n+31)+13\\
\sigma_{\infty}(54n+29)+14\\
\textcolor{blue}{ \sigma_{\infty}(54n+1)+15} \\
\textcolor{red}{ \sigma_{\infty}(54n+41)+16} \\
\sigma_{\infty}(54n+7)+17\\
\sigma_{\infty}(54n+17)+18\\
\sigma_{\infty}(54n+49)+19\\
\vdots \\
\textcolor{black}{{\bf Even_7 \downarrow}} \\
\sigma_{\infty}(54n+11)+1\\
\textcolor{blue}{ \sigma_{\infty}(54n+19)+2} \\
\sigma_{\infty}(54n+23)+3\\
\textcolor{blue}{ \sigma_{\infty}(54n+25)+4} \\
\sigma_{\infty}(54n+53)+5\\
\sigma_{\infty}(54n+13)+6\\
\textcolor{blue}{ \sigma_{\infty}(54n+47)+7} \\
\sigma_{\infty}(54n+37)+8\\
\textcolor{blue}{ \sigma_{\infty}(54n+5)+9} \\
\sigma_{\infty}(54n+43)+10\\
\sigma_{\infty}(54n+35)+11\\
\sigma_{\infty}(54n+31)+12\\
\sigma_{\infty}(54n+29)+13\\
\textcolor{blue}{ \sigma_{\infty}(54n+1)+14} \\
\textcolor{red}{ \sigma_{\infty}(54n+41)+15} \\
\sigma_{\infty}(54n+7)+16\\
\sigma_{\infty}(54n+17)+17\\
\sigma_{\infty}(54n+49)+18\\
\vdots \\
\textcolor{black}{{\bf Odd~d_{7_{next}} \downarrow}} \\
\sigma_{\infty}(54n+11) \\
\textcolor{blue}{ \sigma_{\infty}(54n+19) } \\
\sigma_{\infty}(54n+23) \\
\textcolor{blue}{ \sigma_{\infty}(54n+25) } \\
\sigma_{\infty}(54n+53) \\
\sigma_{\infty}(54n+13) \\
\textcolor{blue}{ \sigma_{\infty}(54n+47) } \\
\sigma_{\infty}(54n+37) \\
\textcolor{blue}{ \sigma_{\infty}(54n+5) } \\
\sigma_{\infty}(54n+43) \\
\sigma_{\infty}(54n+35) \\
\sigma_{\infty}(54n+31) \\
\sigma_{\infty}(54n+29) \\
\textcolor{blue}{ \sigma_{\infty}(54n+1) } \\
\textcolor{red}{ \sigma_{\infty}(54n+41) } \\
\sigma_{\infty}(54n+7) \\
\sigma_{\infty}(54n+17) \\
\sigma_{\infty}(54n+49) \\
\vdots \\
\end{cases}
\begin{cases}
\textcolor{black}{{\bf Odd~d_8\downarrow}} \\
\sigma_{\infty}(54n+17)+2\\
\sigma_{\infty}(54n+49)+3\\
\sigma_{\infty}(54n+11)+4\\
\textcolor{blue}{ \sigma_{\infty}(54n+19)+5} \\
\sigma_{\infty}(54n+23)+6\\
\textcolor{blue}{ \sigma_{\infty}(54n+25)+7} \\
\sigma_{\infty}(54n+53)+8\\
\sigma_{\infty}(54n+13)+9\\
\textcolor{blue}{ \sigma_{\infty}(54n+47)+10} \\
\sigma_{\infty}(54n+37)+11\\
\textcolor{blue}{ \sigma_{\infty}(54n+5)+12} \\
\sigma_{\infty}(54n+43)+13\\
\sigma_{\infty}(54n+35)+14\\
\sigma_{\infty}(54n+31)+15\\
\sigma_{\infty}(54n+29)+16\\
\textcolor{blue}{ \sigma_{\infty}(54n+1)+17} \\
\textcolor{red}{ \sigma_{\infty}(54n+41)+18} \\
\sigma_{\infty}(54n+7)+19\\
\vdots \\
\textcolor{black}{{\bf Even_8 \downarrow}} \\
\sigma_{\infty}(54n+17)+1\\
\sigma_{\infty}(54n+49)+2\\
\sigma_{\infty}(54n+11)+3\\
\textcolor{blue}{ \sigma_{\infty}(54n+19)+4} \\
\sigma_{\infty}(54n+23)+5\\
\textcolor{blue}{ \sigma_{\infty}(54n+25)+6} \\
\sigma_{\infty}(54n+53)+7\\
\sigma_{\infty}(54n+13)+8\\
\textcolor{blue}{ \sigma_{\infty}(54n+47)+9} \\
\sigma_{\infty}(54n+37)+10\\
\textcolor{blue}{ \sigma_{\infty}(54n+5)+11} \\
\sigma_{\infty}(54n+43)+12\\
\sigma_{\infty}(54n+35)+13\\
\sigma_{\infty}(54n+31)+14\\
\sigma_{\infty}(54n+29)+15\\
\textcolor{blue}{ \sigma_{\infty}(54n+1)+16} \\
\textcolor{red}{ \sigma_{\infty}(54n+41)+17} \\
\sigma_{\infty}(54n+7)+18\\
\vdots \\
\textcolor{black}{{\bf Odd~d_{8_{next}} \downarrow}} \\
\sigma_{\infty}(54n+17) \\
\sigma_{\infty}(54n+49) \\
\sigma_{\infty}(54n+11) \\
\textcolor{blue}{ \sigma_{\infty}(54n+19) } \\
\sigma_{\infty}(54n+23) \\
\textcolor{blue}{ \sigma_{\infty}(54n+25) } \\
\sigma_{\infty}(54n+53) \\
\sigma_{\infty}(54n+13) \\
\textcolor{blue}{ \sigma_{\infty}(54n+47) } \\
\sigma_{\infty}(54n+37) \\
\textcolor{blue}{ \sigma_{\infty}(54n+5) } \\
\sigma_{\infty}(54n+43) \\
\sigma_{\infty}(54n+35) \\
\sigma_{\infty}(54n+31) \\
\sigma_{\infty}(54n+29) \\
\textcolor{blue}{ \sigma_{\infty}(54n+1) } \\
\textcolor{red}{ \sigma_{\infty}(54n+41) } \\
\sigma_{\infty}(54n+7) \\
\vdots \\
\end{cases}
\begin{cases}
\textcolor{black}{{\bf Odd~d_9\downarrow}} \\
\textcolor{red}{ \sigma_{\infty}(54n+41)+2} \\
\sigma_{\infty}(54n+7)+3\\
\sigma_{\infty}(54n+17)+4\\
\sigma_{\infty}(54n+49)+5\\
\sigma_{\infty}(54n+11)+6\\
\textcolor{blue}{ \sigma_{\infty}(54n+19)+7} \\
\sigma_{\infty}(54n+23)+8\\
\textcolor{blue}{ \sigma_{\infty}(54n+25)+9} \\
\sigma_{\infty}(54n+53)+10\\
\sigma_{\infty}(54n+13)+11\\
\textcolor{blue}{ \sigma_{\infty}(54n+47)+12} \\
\sigma_{\infty}(54n+37)+13\\
\textcolor{blue}{ \sigma_{\infty}(54n+5)+14} \\
\sigma_{\infty}(54n+43)+15\\
\sigma_{\infty}(54n+35)+16\\
\sigma_{\infty}(54n+31)+17\\
\sigma_{\infty}(54n+29)+18\\
\textcolor{blue}{ \sigma_{\infty}(54n+1)+19} \\
\vdots \\
\textcolor{black}{{\bf Even_9 \downarrow}} \\
\textcolor{red}{ \sigma_{\infty}(54n+41)+1} \\
\sigma_{\infty}(54n+7)+2\\
\sigma_{\infty}(54n+17)+3\\
\sigma_{\infty}(54n+49)+4\\
\sigma_{\infty}(54n+11)+5\\
\textcolor{blue}{ \sigma_{\infty}(54n+19)+6} \\
\sigma_{\infty}(54n+23)+7\\
\textcolor{blue}{ \sigma_{\infty}(54n+25)+8} \\
\sigma_{\infty}(54n+53)+9\\
\sigma_{\infty}(54n+13)+10\\
\textcolor{blue}{ \sigma_{\infty}(54n+47)+11} \\
\sigma_{\infty}(54n+37)+12\\
\textcolor{blue}{ \sigma_{\infty}(54n+5)+13} \\
\sigma_{\infty}(54n+43)+14\\
\sigma_{\infty}(54n+35)+15\\
\sigma_{\infty}(54n+31)+16\\
\sigma_{\infty}(54n+29)+17\\
\textcolor{blue}{ \sigma_{\infty}(54n+1)+18} \\
\vdots \\
\textcolor{black}{{\bf Odd~d_{9_{next}} \downarrow}} \\
\textcolor{red}{ \sigma_{\infty}(54n+41) } \\
\sigma_{\infty}(54n+7) \\
\sigma_{\infty}(54n+17) \\
\sigma_{\infty}(54n+49) \\
\sigma_{\infty}(54n+11) \\
\textcolor{blue}{ \sigma_{\infty}(54n+19) } \\
\sigma_{\infty}(54n+23) \\
\textcolor{blue}{ \sigma_{\infty}(54n+25) } \\
\sigma_{\infty}(54n+53) \\
\sigma_{\infty}(54n+13) \\
\textcolor{blue}{ \sigma_{\infty}(54n+47) } \\
\sigma_{\infty}(54n+37) \\
\textcolor{blue}{ \sigma_{\infty}(54n+5) } \\
\sigma_{\infty}(54n+43) \\
\sigma_{\infty}(54n+35) \\
\sigma_{\infty}(54n+31) \\
\sigma_{\infty}(54n+29) \\
\textcolor{blue}{ \sigma_{\infty}(54n+1) } \\
\vdots \\
\end{cases}
\end{equation}
}
\footnote{All rights reserved. $\copyright$ Michael A. Idowu, 2014.}
\end{landscape}

\begin{conjecture} \label{conj3}
An irrefutable proof of the Collatz conjecture essentially requires 18 fundamental formulae that represent the total stopping time functions of the fundamental covering system of odd integers.
\end{conjecture}

\section{Conclusions}\label{Concl}
A proposed proof of the CC essentially requires deriving the formulae for fundamental total stopping time functions for all odd integers or proving that all the $[3]_{36n+4x}$ numbers eventually converge below these start points using the proposed Collatz based number system, which is both visually demonstrable and theoretically evident.

The proposed covering system of the generalised Collatz based number system requires about 162 distinct sets of odd numbers, from which any other integers could be derived. 
Each fundamental set of odd numbers corresponds to a single fundamental total stopping time function in the proposed schemata.

An irrefutable proof of the Collatz conjecture essentially requires 18 fundamental formulae that represent the total stopping time functions of the fundamental covering system of odd integers.

The Collatz map \ref{CollMachinery1} has many applications. For example, a visual and ingenous method to classify odd numbers to appropriate residue classes modulo 18 is easy. For example, $\{ 349525,1 \} \in [1]_{18}$: $3+4+9+5+2+5 = 28 \equiv 2+8 \equiv 1+0 \in [1]_{18}$; $\{341, 17\} \in [17]_{18}$: $341 \equiv 3+4+1 \equiv  1+7 \in [17]_{18}$. The reader is encouraged to try out this simple technique.

This foundational paper may be regarded as a proposed ``new mathematics'' of the Collatz based number system.  

The whole idea may be used as a new theoretical framework for teaching and understanding elementary number system.

This novel theoretical framework is anticipated to open up new research and further development opportunities in number theory, dynamical systems, and discrete mathematics, including deterministic modelling, metamathematical and optimised integer factorisation, ergodic theory, dynamical systems, covering systems, cryptosystems and cryptography.

One of our aspirations is to further exploit and innovate the main results in visualisable algorithm development.

$\newline$
Dr. Michael A. Idowu
$\newline$
Researcher
$\newline$
Complex Systems Modelling
$\newline$
Abertay University, Dundee DD1 1HG, UK. 
$\newline$
michade@hotmail.com;~m.idowu@abertay.ac.uk


\begin{thebibliography}{99}

\bibitem{Ter76}
Terras, R.; \textit{A stopping time problem on the positive integers}, Acta Arithmetica 30, 241-252 (eng), 1976. http://eudml.org/doc/205476

\bibitem{Lag85} 
 Lagarias, Jeffrey C.; \textit{The 3x + 1 problem and its generalization}, Amer. Math. Monthly 92 (1985), 3-23.

\bibitem{And98}
Andrei, Stefan; Masalagiu, Cristian; \textit{About the Collatz conjecture}. Acta Informatica 35 (2): 167, 1998. doi:10.1007/s002360050117

\bibitem{Van05}
Bendegem, Van; Paul, Jean;  \textit{The Collatz Conjecture: A Case Study in Mathematical Problem Solving}, Logic and Logical Philosophy, volume 14, 7-23, 2005.

\bibitem{Marc96}
Chamberland, Marc;  \textit{A continuous extension of the 3x + 1 problem to the real line}. Dynam. Contin. Discrete Impuls Systems 2: 4, 495-509, 1996.

\bibitem{Gar81}
Garner, Lynn E;  \textit{On the Collatz 3n + 1 Algorithm}. Proceedings of the American Mathematical Society 82 (1): 19-22. doi:10.2307/2044308. JSTOR 2044308, 1981.

\bibitem{Sim99}
Letherman, Simon; Schleicher, Dierk; Wood, Reg;  \textit{The (3n+1)-Problem and Holomorphic Dynamics}. Experimental Mathematics 8: 3, 241-252, 1999.

\bibitem{Mad97}
Maddux, Cleborne D.; Johnson, D. Lamont; ``Logo: A Retrospective''. New York: Haworth Press. p. 160. ISBN 0-7890-0374-0, 1997.

\bibitem{Syr01}
 Lagarias, Jeffrey C.; ``Syracuse problem", in Hazewinkel, Michiel, Encyclopedia of Mathematics, Springer, ISBN 978-1-55608-010-4, 2001.

\bibitem{Ste77}
Steiner, R. P.;  \textit{A theorem on the syracuse problem}, Proceedings of the 7th Manitoba Conference on Numerical Mathematics, pages 553-559, 1977.

\bibitem{Bel06}
Belaga, Edward G.; Mignotte, Maurice;  \textit{Walking Cautiously into the Collatz Wilderness: Algorithmically, Number Theoretically, Randomly}, Fourth Colloquium on Mathematics and Computer Science : Algorithms, Trees, Combinatorics and Probabilities, September 18–22, Institut Élie Cartan, Nancy, France, 2006.

\bibitem{Bel98}
Belaga, Edward G.; Mignotte, Maurice;  \textit{Embedding the 3x+1 Conjecture in a 3x+d Context}, Experimental Mathematics, volume 7, issue 2, 1998.


\bibitem{Sim05}
Simons, J.; de Weger, B.;  \textit{Theoretical and computational bounds for m-cycles of the 3n + 1 problem}, Acta Arithmetica (on-line version 1.0, November 18, 2003), 2005.

\bibitem{Sin10}
Sinyor, J.;, \textit{"The 3x+1 Problem as a String Rewriting System}, International Journal of Mathematics and Mathematical Sciences, volume 2010, Article ID 458563, 6 pages, 2010.


\bibitem{Lag11} 
 Lagarias, Jeffrey C.; \textit{The 3x + 1 problem: An annotated bibliography (19631999) (sorted by author)}, Arxiv (1985), 2003. http://arxiv.org/abs/math/0309224.

\bibitem{Lag12}
 Lagarias, Jeffrey C.; \textit{The 3x+1 problem: An annotated bibliography, II (2000-2009)}. Arxiv. Available at arxiv.org/abs/math/0608208.
 
\end{thebibliography}
\end{document}